\documentclass[10pt]{article}
\usepackage{amsmath,amssymb}
\newcommand{\bra}{\langle}
\newcommand{\ket}{\rangle}
\newcommand{\ep}{\hfill {$\Box$}}

\newcommand{\pf }{{\noindent{\bf Proof. }}}

\newtheorem{thm}{Theorem}[section]
\newtheorem{cor}{Corollary}[section]
\newtheorem{defin}{Definition}[section]
\newtheorem{lem}{Lemma}[section]
\newtheorem{prop}{Proposition}[section]

\title{Spectral stability for compact perturbations of Toeplitz matrices}
\author{\small{Astaburuaga M.A., Bourget O., Cort\'es V.H.\footnote{Supported by the Grants Fondecyt 1120786 and Conicyt PIA-ACT1112}}\\
\small{Facultad de Matem\'aticas, Pontificia Universidad Cat\'olica de Chile,}\\
\small{Av. Vicu\~na Mackenna 4860, Macul, Santiago, Chile}\\
\small{bourget@mat.puc.cl}\\
\small{phone: (56 2) 2354 4509}\\
\small{fax: (56 2) 2552 5916}
}
\date{}
\begin{document}
\maketitle
%\newpage

\begin{abstract}
Let $f$ be a regular real-valued non-constant symbol defined on the one dimensional torus ${\mathbb T}$. Denote respectively by $\kappa$ and $T$, its set of critical points and the associated Toeplitz matrix on $l^2({\mathbb N})$. If $V$ is a suitable compact perturbation, we prove that the operator $T+V$ has no singular continuous spectrum and only finite point spectrum away from the set of thresholds $f(\kappa)$. We also obtain some propagation estimates and apply these results to concrete examples.
\end{abstract}


%**************************************
\section{Introduction}

Toeplitz operators have been studied extensively over the years and their spectral properties are well-known, see e.g. \cite{bg}, \cite{bs}, \cite{dp}, \cite{gs}, \cite{rr}. Their compact perturbations have also been characterized \cite{xia}. However, aside from the invariance of the essential spectrum, few information is available concerning the stability of these spectral properties under compact perturbations, see e.g. \cite{john}.

In this paper, we show that under suitable regularity conditions on the symbol and the perturbation, the spectral properties of the associated Toeplitz operator remain qualitatively stable away from a set of critical values (see Theorem \ref{toeplitz}). This result is based on a positive commutator technique (regular Mourre theory). We apply it to various contexts. First, we consider compact perturbations of the discrete Schr\"odinger operator on the half-line and deduce Theorem \ref{schro} as a counterpart of the results obtained in \cite{sah} on the lattice. Second, we consider the case of suitable finite rank perturbations and extend the results obtained in \cite{john} for the rank one case (Theorem \ref{rank1}). Finally, we study the spectral properties of symmetrized products of Toeplitz operators in Theorem \ref{product}.

These results are introduced in Section 2. We dedicate Section 3 to their proofs and provide a dynamical interpretation in Section 4 with some propagation estimates. All these results are contrasted with the case of Laurent operators defined on the lattices ${\mathbb Z}$ and ${\mathbb Z}^d$ in Section 5.

\noindent{\bf Notations.} If ${\cal H}$ denotes an infinite-dimensional (complex) Hilbert space, ${\cal B}({\cal H})$ is the algebra of bounded operators acting on ${\cal H}$. The resolvent set of an operator $B$ in ${\cal B}({\cal H})$ is denoted by $\rho(B)$ and its spectrum by $\sigma(B):= {\mathbb C}\setminus \rho(B)$. The one-dimensional torus is denoted by ${\mathbb T}:={\mathbb R}/2\pi {\mathbb Z}$. If $H$ is a self-adjoint operator defined on ${\cal H}$, its spectral family is denoted by $(E_{\Delta}(H))_{\Delta \in {\cal B}({\mathbb R})}$, where ${\cal B}({\mathbb R})$ stands for the family of Borel sets of ${\mathbb R}$. 
The continuous and point subspaces of the operator $H$ are respectively denoted by ${\cal H}_c(H)$ and ${\cal H}_{pp}(H)$. In this paper, we deal mostly with the Hilbert spaces
\begin{eqnarray*}
l^{2} ({\mathbb N}) &=& \{ \psi = (\psi_n)_{n\in {\mathbb N}} : \sum_{n\in {\mathbb N}} |\psi_n |^2  < \infty  \}\\
\mbox{and}\quad l^{2} ({\mathbb Z}^d) &=& \{ \psi = (\psi_n)_{n\in {\mathbb Z}^d} : \sum_{n\in {\mathbb Z}^d} |\psi_n |^2  < \infty  \} \quad\text{for}\, d\in {\mathbb N}
\end{eqnarray*}
equipped with their usual inner product $\bra \varphi,\psi \ket_{l^{2} ({\mathbb N})} = \sum_{n\in {\mathbb N}} \bar{\varphi}_n \psi_n$ and $\bra \varphi,\psi \ket_{l^{2} ({\mathbb Z}^d)} = \sum_{n\in {\mathbb Z}^d} \bar{\varphi}_n \psi_n$. The canonical orthonormal bases of $l^{2} ({\mathbb N})$ and $l^2({\mathbb Z}^d)$ are denoted by $(e_n)_{n\in {\mathbb N}}$ and $(e_n)_{n\in {\mathbb Z}^d}$. If ${\cal N}$ is a subset of ${\mathbb N}$ or ${\mathbb Z}^d$, the linear span of the vectors $(e_n)_{n\in {\cal N}}$ is denoted by $\bra e_n; n\in {\cal N}\ket$ and its Hilbert closure by $\overline{\bra e_n; n\in {\cal N}\ket}$.

$L^{2} ({\mathbb T}^d)$ denotes the Hilbert space of square integrable complex functions with inner product
\begin{equation*}
\langle f, g  \rangle =  \frac{1}{(2\pi)^d  } \,  \int_{{\mathbb T}^d} \overline{f(\theta)} g(\theta) \, d\theta \, ,
\end{equation*}
With these notations, the Fourier transform  ${\mathcal F} : L^2({\mathbb T}^d) \to l^{2} ({\mathbb Z}^d)$ is defined by: ${\mathcal F} f = (\hat{f}_n)_{n \in {\mathbb Z}^d}$ where $  \hat{f}_n$ is the $n$-th Fourier coefficient of the function $f\in L^{2} ({\mathbb T}^d)$:
\begin{equation}\label{fourier1}
\hat{f}_n = \frac{1}{(2 \pi)^d} \int_{{\mathbb T}^d} e^{-in \cdot \theta} f(\theta) \, d\theta \, .
\end{equation}
The spaces $L^{\infty} ({\mathbb T}^d)$, $C^0({\mathbb T}^d)$ and $C^k({\mathbb T}^d)$ ($k\in {\mathbb N}$) stand respectively for the linear spaces of essentially bounded complex functions, continuous complex functions and $k$-th continuously differentiable complex functions defined on ${\mathbb T}^d$. The Wiener algebra is denoted by: ${\cal A}({\mathbb T}^d):= \{f\in L^{\infty}({\mathbb T}^d); (\hat{f}_n)_{n \in {\mathbb Z}^d}\in l^1({\mathbb Z}^d)\}$.

%*************************
\section{Main results}

The main ingredients are introduced in Sections 2.1 and 2.2. The results are stated in Sections 2.3 and 2.4.

%**********************************************
\subsection{Laurent, Toeplitz and Hankel operators}

\paragraph{Laurent operators.} Let $f\in L^{\infty}({\mathbb T})$. The Laurent operator associated to $f$ is defined by: $L_f: l^{2} ({\mathbb Z}) \rightarrow l^{2} ({\mathbb Z}),$
\begin{equation}\label{laurentope1}
(L_f \varphi)_n=  (({\cal F}f)*\varphi)_n = \sum_{k \in {\mathbb Z}} \hat{f}_{n-k}\varphi_k  \, ,
\end{equation}
for all $\varphi\in \bra e_n;n\in {\mathbb Z}\ket$ and all $n\in {\mathbb Z}$. Since $f\in L^{\infty}({\mathbb T})$, $L_f$ extends as a bounded operator on $l^{2} ({\mathbb Z})$ and $\|L_f \|=\|f\|_{\infty}$. The operator $L_f$ is unitarily equivalent to the multiplication operator by the function $f$ on $L^2({\mathbb T})$: $L_f={\mathcal F}f{\mathcal F}^*$. For any functions $f$ and $g$ in $L^{\infty}({\mathbb T})$ and any $c\in {\mathbb C}$, we have that: $L_{f+g}=L_f+L_g$, $L_{fg}=L_f L_g$, $L_{cf}=cL_f$, $L_f^*=L_{\bar{f}}$. In particular, $[L_f,L_g]=0$. $L_1$ is the identity on $l^{2} ({\mathbb Z})$.

\paragraph{Toeplitz operators.} Let $f\in L^{\infty}({\mathbb T})$. The Toeplitz operator associated to $f$ is defined by: $T_f: l^{2} ({\mathbb N}) \rightarrow l^{2} ({\mathbb N}),$
\begin{equation}\label{laurentope1}
(T_f \varphi)_n= \sum_{k \in {\mathbb N}} \hat{f}_{n-k}\varphi_k  \, ,
\end{equation}
for all $\varphi\in \bra e_n; n\in {\mathbb N}\ket$ and all $n\in {\mathbb N}$. Since $f\in L^{\infty}({\mathbb T})$, $T_f$ extends as a bounded operator on $l^{2} ({\mathbb N})$, $\|T_f\|\leq \|f\|_{\infty}$. For any functions $f$ and $g$ in $L^{\infty}({\mathbb T})$ and any $c\in {\mathbb C}$, we have that: $T_{f+g}=T_f+T_g$, $T_{cf}=cT_f$, $T_f^*=T_{\bar{f}}$. But $T_{fg}\neq T_f T_g$ and $[T_f,T_g]\neq 0$ in general (see e.g. formula (\ref{sarason})). $T_1$ is the identity on $l^{2} ({\mathbb N})$. If the function $f\in L^{\infty}({\mathbb T})$ is real-valued, then $T_f$ is self-adjoint, $\sigma(T_f)$ is equal to the essential range of $f$ \cite{hw}, \cite{wid} and the operator $T_f$ has purely absolutely continuous spectrum \cite{ros}.

With the direct sum $l^2({\mathbb Z})= \overline{\bra e_n; n\geq 1\ket} \oplus \overline{\bra e_n; n\leq 0 \ket}$, we observe that $l^2({\mathbb N})$ and the closed subspace $\overline{\bra e_n; n\geq 1\ket}$ of $l^2({\mathbb Z})$ are canonically unitarily isomorphic, which provides a natural embedding of $l^2({\mathbb N})$ into $l^2({\mathbb Z})$. If $P$ denotes the orthogonal projection on $\overline{\bra e_n; n\geq 1\ket}$ in $l^2({\mathbb Z})$ and $P^{\perp}=I-P$ the orthogonal projector on its orthocomplement, we identify the operators $T_f$ and $PL_fP$. Note that $l^2({\mathbb N})$ and the closed subspace $\overline{\bra e_n; n\leq 0 \ket}$ of $l^2({\mathbb Z})$ can also be identified naturally.

\paragraph{Hankel operators.} Let $f\in L^{\infty}({\mathbb T})$. We define the Hankel operator associated to $f$ by: $H_f: l^{2} ({\mathbb N}) \rightarrow l^{2} ({\mathbb N}),$
\begin{equation*}
(H_f \varphi)_n= \sum_{k \in {\mathbb N}} \hat{f}_{n+k-1}\varphi_k  \, ,
\end{equation*}
for all $\varphi\in \bra e_n;n\in {\mathbb N}\ket$ and all $n\in {\mathbb N}$. $H_f$ extends also as a bounded operator on $l^{2} ({\mathbb N})$ (Nehari Theorem \cite{neh}). If in addition $f\in C^0({\mathbb T})$, $H_f$ is compact (Hartman Theorem \cite{bp}, \cite{hart}). For more details see also \cite{part}, \cite{pell}, \cite{power}.

Once identified the closed subspaces $\overline{\bra e_n; n\geq 1\ket}$ and $\overline{\bra e_n; n\leq 0\ket}$ of $l^2({\mathbb Z})$ with copies of $l^2({\mathbb N})$ as stated above, we also identify the operator $P^{\perp} L_f P$ (resp. $P L_f P^{\perp}$) and the Hankel operator $H_f$ (resp. $H_{\bar{f}}^*$). All these correspondences are used freely in this paper: namely, we will write for all $f\in L^{\infty}({\mathbb T}),$
\begin{equation}
T_f = PL_fP\, , \quad H_f = P^{\perp} L_f P\, , \quad H_{\bar{f}}^* =PL_f P^{\perp} \label{identif}
\end{equation}
We deduce in particular that for any $(f, g) \in L^{\infty}(\mathbb{T})\times L^{\infty}(\mathbb{T})$, one has that
\begin{equation}\label{sarason}
T_f T_g = PL_fPL_gP=PL_{fg}P - PL_f P^{\perp}L_gP= T_{fg}- H_{\bar{f}}^* H_g \, .
\end{equation}

\paragraph{Symbols.} The function $f\in L^{\infty}({\mathbb T})$ associated to the operators $L_f$, $T_f$ and $H_f$ above is sometimes called the symbol of these operators. For real-valued functions $f$ in $L^{\infty}({\mathbb T})$, $\bar{\hat{f}}_n = \hat{f}_{-n}$ and the operators $L_f$ and $T_f$ are self-adjoint. The set of critical points of $f$ is denoted by
\begin{equation*}
\kappa_f = \{\theta \in {\mathbb T}; \text{f is not differentiable at } \theta \text{ or } f'(\theta)=0\} \, .
\end{equation*}
In particular, for $f\in C^1({\mathbb T})$ real-valued, the set $\kappa_f$ is a compact subset of ${\mathbb T}$ and the set of thresholds $f(\kappa_f)$ is a compact subset of $\sigma (T_f)=\mathrm{Ran} f$.

%*************************************
\subsection{Commutation and regularity}

Let ${\cal H}$ be a Hilbert space. Consider $A$ a self-adjoint operator defined on $\mathcal H$ with domain ${\mathcal D}(A)$ and let $B$ be a bounded operator on $\mathcal H$. We say that $B$ is of class $C^1$ w.r.t. $A$ if the sesquilinear form $F$ defined by
\[F(\varphi,\psi):=\langle A\varphi,B\psi\rangle - \langle \varphi,BA\psi\rangle\]
for any $(\varphi,\psi)\in {\cal D}(A)\times {\cal D}(A)$, extends continuously to a bounded form on ${\cal H}\times {\cal H}$ (w.r.t. the product topology of ${\cal H}\times {\cal H}$). The (unique) bounded linear operator associated to the extension is denoted by $\mathrm{ad}_A (B) =[A,B]$. We denote: $C^1(A):= \{B\in {\cal B}({\cal H}); B$ is of class $C^1$ w.r.t $A \}$.
In practice, it is enough to check this continuity property on some core of $A$.

Higher order commutators are defined inductively as follows. With the convention that $C^0(A)={\cal B}({\cal H})$ and $\mathrm{ad}^0_A B=B$ for all $B\in {\cal B}({\cal H})$, we say that the operator $B \in {\mathcal B}(\mathcal H)$ is of class $C^k$ w.r.t. $A$ for some $k\in {\mathbb N}$, if $B\in C^{k-1}(A)$ and $\text{ad}_A^{k-1} B \in C^1(A)$. We denote: $\mathrm{ad}_A (\mathrm{ad}_A^{k-1}(B)) =\mathrm{ad}_A^{k}(B)$ and $C^k(A):= \{B\in {\cal B}({\cal H}); B$ is of class $C^k$ w.r.t $A \}$. Also, $C^{\infty}(A):= \cap_{k\in {\mathbb N}}C^k(A)$.

\noindent{\bf Remark:} We can also consider fractional order regularities \cite{abmg}, \cite{sah}. Let us just mention that for $B\in {\cal B}({\cal H})$, we say that:
\begin{itemize}
\item $B\in {\cal C}^{0,1}(A)$ if:
\begin{equation*}
\int_0^1 \|e^{iA\tau}B e^{-iA\tau}-B\|\,\frac{d\tau}{\tau} < \infty \, .
\end{equation*}
\item $B\in {\cal C}^{1,1}(A)$ if:
\begin{equation*}
\int_0^1 \|e^{iA\tau}B e^{-iA\tau}+e^{-iA\tau}B e^{iA\tau} -2B\|\,\frac{d\tau}{\tau^2} < \infty \, .
\end{equation*}
\end{itemize}
Clearly, ${\cal C}^{0,1}(A)$ and ${\cal C}^{1,1}(A)$ are linear subspaces of ${\cal B}({\cal H})$. They are stable under adjunction $*$. It is also known that if $B\in C^1(A)$ and $\mathrm{ad}_A B\in {\cal C}^{0,1}(A)$, then $B\in {\cal C}^{1,1}(A)$ and that $C^2(A)\subset {\mathcal C}^{1,1}(A)\subset C^1(A)$ (see e.g. inclusions 5.2.19 in \cite{abmg}).

Let $B$ be a self-adjoint operator defined on the Hilbert space ${\cal H}$. We will say that a Limiting Absorption Principle (LAP) holds for $B$ on some Borel subset $\Lambda\subset {\mathbb R}$ (w.r.t some auxiliary self-adjoint operator $A$) if statements (a)--(c) below are satisfied:
\begin{itemize}
\item[(a)] For any compact subset $K\subset \Lambda$
\begin{equation*}
\sup_{ \Im z \neq 0, \Re z \in K} \|\bra A \ket^{-1}(z-B)^{-1}\bra A \ket^{-1} \| < \infty \, .
\end{equation*}
\item[(b)] If $z$ tends to $\lambda \in \Lambda$ (non-tangentially), then $\bra A \ket^{-1}(z-B)^{-1}\bra A \ket^{-1}$ converges in norm to a bounded operator denoted $F^+(\lambda)$ (resp. $F^-(\lambda)$) if $\Im z >0$ (resp. $\Im z <0$). This convergence is uniform on any compact subset $K\subset \Lambda$.
\item[(c)] The operator-valued functions defined by $F^{\pm}$ are continuous on each connected component of $\Lambda$, w.r.t. the norm topology on ${\cal B}({\cal H})$.
\end{itemize}

\noindent{\bf Remark:} If $\Lambda\subset \rho(B)$, the properties described above are trivially satisfied. In this case, $F^+(\lambda)=F^-(\lambda)=\bra A \ket^{-1}(\lambda-B)^{-1}\bra A \ket^{-1}$ for any $\lambda\in \Lambda$.

Note that for suitable self-adjoint operators $B$, Mourre Theory relates the existence of commutation and conjugacy properties (see Section 3.1) to the existence of a LAP. The choice of the operator $A$ depends in practice on $B$. In the next section, we specify such a choice for Toeplitz operators.

%*****************************
\subsection{Conjugate operators}

Let ${\mathfrak X}$ denote the linear operator defined on the canonical orthonormal basis of $l^2({\mathbb Z})$ by: ${\mathfrak X} e_n = n e_n$, $n\in {\mathbb Z}$. It is essentially self-adjoint on $\bra e_n; n\in {\mathbb Z}\ket$ and its self-adjoint extension with domain ${\mathcal D}_{\mathfrak X}= \{\psi=(\psi_n)_{n\in {\mathbb Z}}   \in l^2({\mathbb Z}) :  \sum_{ n\in {\mathbb Z}} n^2 |\psi_n|^2 < \infty \}$ is also denoted ${\mathfrak X}$. Similarly, let $X$ denote the linear operator defined on the canonical orthonormal basis of $l^2({\mathbb N})$ by: $X e_n = n e_n$, $n\in {\mathbb N}$. It is essentially self-adjoint on $\bra e_n; n\in {\mathbb N}\ket$ and its self-adjoint extension with domain ${\mathcal D}_{X}= \{\psi=(\psi_n)_{n\in {\mathbb N}} \in l^2({\mathbb N}) :  \sum_{ n\in {\mathbb N}} n^2 |\psi_n|^2 < \infty \}$ is also denoted $X$. The operators ${\mathfrak X}$ and $X$ are respectively the position operators on $l^2({\mathbb Z})$ and $l^2({\mathbb N})$. We observe that $P\in C^1({\mathfrak X})$ and that $[{\mathfrak X},P]=0$. So, ${\mathfrak X}= P {\mathfrak X} P \oplus P^{\perp} {\mathfrak X} P^{\perp}$. In the following, we identify the operators $X$ and $P{\mathfrak X}P$.

Note that ${\cal F}^*{\mathfrak X} {\cal F}= -i\partial_{\theta}$. By Fourier transform (\ref{fourier1}), we deduce that for $h\in C^k({\mathbb T})$, $k\in {\mathbb N}$, $L_h\in C^k({\mathfrak X})$ and for all $j\in \{0,\ldots,k\}$, $\mathrm{ad}_{\mathfrak X}^j L_h = (-i)^j L_{h^{(j)}}$. Using the identifications $X=P{\mathfrak X}P$ and $T=PLP$, we deduce that:
\begin{lem}\label{ckx} Let $h \in C^k({\mathbb T})$. Then, $T_h \in C^k(X)$ and for all $j\in \{0,\ldots,k\}$, $\mathrm{ad}_X^j T_h= (-i)^j T_{h^{(j)}}$.
\end{lem}

Let $g$ be a real-valued function which belongs to $C^2({\mathbb T})$. By Lemma \ref{ckx}, the corresponding Toeplitz operator $T_g$ belongs to $C^1(X)$, so $T_g {\cal D}_X \subset {\cal D}_X$. This allows us to define on ${\cal D}_X$ the symmetric operator $A_g$ by:
\begin{equation}\label{Ag}
A_g =   \frac{1}{2}(T_g X + XT_g)\, .
\end{equation}
According to Lemma \ref{ckx}, we can rewrite:
\begin{equation}\label{Ag2}
A_g = T_g X + \frac{1}{2}\mathrm{ad}_X T_g = XT_g - \frac{1}{2}\mathrm{ad}_X T_g \, .
\end{equation}

\begin{lem}\label{sacom} Let $g$ in $C^2({\mathbb T})$ be real-valued. Then, the symmetric operator $A_g$ is essentially self-adjoint. The linear space $\bra e_n; n\in {\mathbb N}\ket$ is a core for $A_g^{**}$.
\end{lem}
\pf The linear space ${\cal S}:=\bra e_n; n\in {\mathbb N}\ket$ is a core for $X$ and $X^2$. The Toeplitz operator $T_g$ is bounded symmetric hence self-adjoint on $l^2({\mathbb N})$. It implies that the operator $A_g$ and its restriction to ${\cal S}$ (denoted by $A_g^o$) are also symmetric. We also notice that for all $\varphi \in {\cal S}$, $\|\varphi \| \leq \|X\varphi \| \leq \|X^2\varphi \|$. From (\ref{Ag2}) and Lemma \ref{ckx}, we deduce that for some $C>0$ and for any $\varphi \in {\cal S}$, $\|A_g^o \varphi \| \leq C \|X^2\varphi \|$. Using again Lemma \ref{ckx}, the following identity holds as a sesquilinear form on ${\cal S}\times {\cal S}$:
\begin{eqnarray*}
[A_g^o,X^2] &=& -\frac{1}{2}\left(X^2[X,T_g] +2X[X,T_g]X +[X,T_g]X^2\right)\\
&=& \frac{1}{2}\left([X,[X,T_g]]X -4X[X,T_g]X -X[X,[X,T_g]]\right)= \frac{1}{2}\left(-T_{g''}X +4i XT_{g'}X +XT_{g''}\right)\, .
\end{eqnarray*}
With the same arguments, we deduce that $|\bra A_g^o \varphi, X \varphi  \ket -  \bra  X \varphi,  A_g^o\varphi  \ket| \leq C \|  X\varphi \|^2$ for some $C>0$ and all $\varphi \in {\cal S}$. By Theorem X.37 in \cite{rs2}, $A_g^o$ is essentially self-adjoint. $A_g$ is a symmetric extension of $A_g^o$. This implies that $A_g$ is also essentially self-adjoint and that $A_g^{**}=A_g^{o\,**}$ (see e.g. \cite{rs2} Section VIII.2). \ep

In the remainder of the paper, we abuse notations and denote by $A_g$ the  self-adjoint extension $A_g^{**}$. We conclude this section with the following observation:
\begin{lem}\label{relbound} Let $g$ in $C^2({\mathbb T})$ be real-valued. Then, the operators
$A_g X^{-1}$ and $A_g^2 X^{-2}$ are bounded.
\end{lem}

%***************************
\subsection{An abstract result}

Let $f\in C^0({\mathbb T})$ be real-valued. We know that $\sigma(T_f)=\sigma_{ess}(T_f)=$ Ran $f$. It follows immediately from Weyl Theorem that for any compact symmetric operator $V$, $\sigma_{ess}(T_f+V)=\sigma_{ess}(T_f)=$ Ran $f$. In addition, we have that:
\begin{thm}\label{toeplitz} Consider a non-constant real-valued symbol $f\in C^3({\mathbb T})$. Let $H=T_f+V$ with $V$ a compact symmetric operator defined on $l^2({\mathbb N})$ such that
$V\in {\cal C}^{1,1}(A_{f'})$. Then,
\begin{itemize}
\item[(a)] given any Borel set $\Lambda$ such that $\overline{\Lambda}\subset$ Ran $f\setminus f(\kappa_f)$, $H$ has at most a finite number of eigenvalues in $\Lambda$. Each of these eigenvalues has finite multiplicity.
\item[(b)] a LAP holds for $H$ on Ran $f\setminus \sigma_{\text{pp}}(H)\cup f(\kappa_f)$ w.r.t $A_{f'}$. $H$ has no singular continuous spectrum in Ran $f\setminus f(\kappa_f)$.
\end{itemize}
In particular, if $f(\kappa_f)$ has a finite number of accumulation points, $H$ has no singular continuous spectrum. 
\end{thm}
See Section 3.4 for the proof. A similar result holds for Laurent operators (see Theorem \ref{laurent-glob}). Note also that:
\begin{itemize}
\item the regularity hypothesis on $f$ can be relaxed with extra-technicalities.
\item a generalized version of the LAP can also be stated under the hypotheses of Theorem \ref{toeplitz} (see e.g. \cite{abmg} Theorem 7.3.1).
\item stronger regularity assumptions on $f$ and $V$ would imply stronger regularity properties for the spectral measure of $H$ \cite{jmp}.
\item the distribution of the point spectrum is not studied here.
\end{itemize}
We illustrate Theorem \ref{toeplitz} with various applications.

%******************************************************************
\subsection{Examples}
\subsubsection{Discrete Schr\"odinger operator on the half-line}

Consider the function $f\in C^{\infty}({\mathbb T})$ defined by $f(\theta)=2\cos \theta$. The associated Toeplitz operator $T_f$ is defined by $T_f e_1 =e_2$ and $T_f e_n= e_{n-1}+e_{n+1}$ if $n\geq 2$. We have that: $\sigma(T_f)=[-2,2]$ and $\kappa_f=\{0,\pi\}$. So, $f(\kappa_f)=\{-2,2\}$. Let us consider a compact symmetric operator $V$ defined on the canonical orthonormal basis by: $Ve_n =v_n e_n$, $n\geq 1$ where the sequence $(v_n)$ vanishes at infinity. Denoting $H:=T_f+V$, we have that $\sigma_{ess}(H)=[-2,2]$ by Weyl Theorem.

To measure the regularity of the operator $V$ w.r.t. the self-adjoint operator $A_{f'}$ defined by (\ref{Ag}), we introduce the family of norms $(q_k)_{k\geq 0}$ defined on ${\mathbb C}^{\mathbb N}$ by:
\begin{equation}\label{seminorm}
q_0(\gamma) = \| \gamma \|_{\infty}:=  \sup_{n \in {\mathbb N}} |\gamma_n | \quad \text{and}\quad q_{k+1}(\gamma) = q_k(\gamma) + \|\xi^{k+1} \Delta^{k+1} \gamma \|_{\infty}  \, ,
\end{equation}
where $\Delta\gamma$ and $\xi \gamma$ are defined by: $(\Delta \gamma)_n = \gamma_{n}-\gamma_{n+1}$ and $(\xi \gamma)_n =n\gamma_n$, $n\in {\mathbb N}$.
We also introduce the following sets of hypotheses:
\begin{itemize}
\item[{\bf (S)}] There exist $0<a_1<b_1 < \infty$ such that: \begin{equation*}
\int_1^{\infty} \sup_{a_1r\leq n \leq b_1r} |\gamma_n |\, dr < \infty \, .
\end{equation*}
\item[{\bf (M)}] $q_{1}(\gamma)<\infty$, $\lim_{n\rightarrow \infty} \gamma_n=0$ and there exist $0<a_2<b_2 < \infty$ such that: \begin{equation*}
\int_1^{\infty} \sup_{a_2r\leq n \leq b_2r} |\gamma_{n+1} -\gamma_n|\, dr < \infty
\end{equation*}
\item[{\bf (L)}] $q_{2}(\gamma) < \infty$ and $\lim_{n\rightarrow \infty} \gamma_n=0$.
\end{itemize}

\begin{thm}\label{schro} Let $f$ be defined by $f=2\cos$. Let $(s_n)_{n\in {\mathbb N}}$, $(m_n)_{n\in {\mathbb N}}$ and $(l_n)_{n\in {\mathbb N}}$ be real-valued sequences which satisfy conditions {\bf S}, {\bf M} and {\bf L} respectively. Consider $H:=T_f+V$ where $Ve_n =v_n e_n$, $v_n = s_n+m_n+l_n$, $n\in {\mathbb N}$. Then, $\sigma_{ess}(H)=[-2,2]$ and
\begin{itemize}
\item any open interval $\Lambda$ such that $\overline{\Lambda} \subset (-2,2)$ contains at most a finite number of eigenvalues. Each of these eigenvalues has finite multiplicity.
\item a LAP holds for $H$ on $(-2,2) \setminus \sigma_{\text{pp}}(H)$ w.r.t $A_{f'}$. $H$ does not have any singular continuous spectrum.
\end{itemize}
\end{thm}

The proof is developed in Section 3.5. Theorem \ref{schro} is the half-line analog of Theorem 2.1 in \cite{sah} (which is itself a special case of Theorem \ref{laurent-glob}). 

%*****************************************
\subsubsection{Finite rank perturbations}

Our next example is motivated by \cite{john}. We will say that a vector $\psi\in l^2({\mathbb N})$ satisfies the hypothesis {\bf H} if for some $0<a<b<\infty,$
\begin{equation*}
\int_{1}^{\infty} ( \sum_{n\in {\mathbb N}\cap [ar,br]}  | \psi_n |^2 \, )^{1/2}\, dr < \infty \,.
\end{equation*}

\begin{thm}\label{rank1} Let $f \in C^3({\mathbb T})$ be non-constant and real-valued. Let $N\in {\mathbb N}$ and consider a finite family of vectors $(\psi_k)_{k=1}^N \subset l^2({\mathbb N})$ such that for all $k\in \{1,\ldots,N\}$, $\psi_k$ satisfies {\bf H} or belongs to ${\cal D}(A_{f'}^2)$. Let $(H_{\beta})_{\beta \in {\mathbb R}^N}$ be the family of operators defined by $H_{\beta}:=T_f + V_{\beta}$ where
\begin{equation*}
V_{\beta}=\sum_{k=1}^N \beta_k |\psi_k\ket \bra \psi_k |\, .
\end{equation*}
Then for any $\beta \in {\mathbb R}^N$, $\sigma_{ess}(H_{\beta})=$ Ran $f$ and statements (a)-(b) of Theorem \ref{toeplitz} hold for $H_{\beta}$.
\end{thm}
The proof is developed in Section 3.6. Theorem \ref{rank1} extends the rank-one case studied in \cite{john}.

%*****************************************
\subsubsection{On the product of Toeplitz operators}

\begin{thm}\label{product} Let $f,g$ be two real-valued functions in $C^3({\mathbb T})$ such that their product $h=fg$ is not constant. Let $H=\Re (T_fT_g)$. Then, $\sigma_{ess}(H)=$ Ran $h$ and
\begin{itemize}
\item given any Borel set $\Lambda$ such that $\overline{\Lambda}\subset$ Ran $h\setminus h(\kappa_h)$, $H$ has at most a finite number of eigenvalues in $\Lambda$. Each of these eigenvalues has finite multiplicity.
\item a LAP holds for $H$ on Ran $h \setminus \sigma_{\text{pp}}(H) \cup h(\kappa_h)$ w.r.t $A_{h'}$. $H$ has no singular continuous spectrum in Ran $h\setminus h(\kappa_h)$.
\end{itemize}
\end{thm}
From the proof developed in Section 3.7, we could relax the regularity hypotheses on the functions $f$ and $g$. We shall not consider it here. 

%********************
\section{Technicalities}

First, we start by reviewing the main features of the regular Mourre theory for self-adjoint operators that are used in the proof of Theorem \ref{toeplitz}. To grasp an overview of the theory, the reader is referred to Chapter 4 in \cite{cfks} and \cite{gg}. For detailed and optimal results, we refer to Chapter 7 in \cite{abmg}.

We also use the following notations: $S\thickapprox T$ if $S-T$ is a compact operator and $S\lesssim T$ (resp. $S\gtrsim T$) if $S\leq T+K$ (resp. $S\geq T+K$) for some compact operator $K$. For example, given $(f,g)\in L^{\infty}({\mathbb T})\times L^{\infty}({\mathbb T})$, (\ref{sarason}) rewrites $T_fT_g\simeq T_{fg}$ whenever $f$ or $g$ is a continuous function (Hartman Theorem).

%**********************************
\subsection{Regular Mourre Theory}

Throughout this section, ${\cal H}$ is a Hilbert space and $H$ a bounded self-adjoint operator defined on ${\cal H}$. We start by recalling  the concept of conjugacy which is central in our discussion. 

\begin{defin}\label{propagation} Assume that there exist a self-adjoint operator $A$ with domain ${\cal D}(A) \subset {\cal H}$
such that $H\in C^1(A)$. For a given $\Lambda \in {\cal B}({\mathbb R})$, we say that
\begin{itemize}
\item $H$ is weakly conjugate w.r.t. $A$ if $i[A,H] > 0$ i.e. for all $\varphi \in {\cal H}\setminus \{0\}$, $\bra \varphi, i[A,H]\varphi \ket >0$.
\item $H$ is conjugate w.r.t. $A$ on ${\Lambda}$ if there exist $c >0$ such that: $E_{\Lambda}(H) i[A,H] E_{\Lambda}(H) \gtrsim c E_{\Lambda}(H)$.
\item $H$ is strictly conjugate w.r.t. $A$ on ${\Lambda}$ if there exist $c >0$ such that:
$E_{\Lambda}(H) i[A,H] E_{\Lambda}(H) \geq c E_{\Lambda}(H)$.
\end{itemize}
\end{defin}
Mourre Theory (or conjugate operator method) provides a control of the point spectrum via the Virial Theorem, which states that $E_{\{\lambda\}}(H) i[A,H] E_{\{\lambda\}}(H) =0$ for all $\lambda \in {\mathbb R}$ if $H \in C^1(A)$ (see e.g. \cite{gg} or Proposition 7.2.10 in \cite{abmg}). As a consequence, if $H$ is weakly conjugate w.r.t $A$, then $H$ has no eigenvalue. Similarly, if $H$ is strictly conjugate w.r.t $A$ on some Borel set $\Lambda \subset {\mathbb R}$, then $H$ has no eigenvalue in $\Lambda$. The Virial Theorem also implies that:
\begin{prop}\label{virial-2} Let $H\in C^1(A)$. Assume that $H$ is conjugate w.r.t. $A$ on the Borel set $\Lambda\subset {\mathbb R}$.
Then, $H$ has a finite number of eigenvalues in $\Lambda$. Each of these eigenvalues has finite multiplicity.
\end{prop}
See \cite{abmg}, Corollary 7.2.11 for the proof. Mourre Theory provides also the existence of a LAP away from the set of eigenvalues and therefore allows to rule out the existence of singular continuous spectrum. We refer to Section 7.3 in \cite{abmg} for a proof of the following result.
\begin{thm}\label{nosc2} Let $\Lambda\subset {\mathbb R}$ be an open set. Assume that $H\in {\cal C}^{1,1}(A)$ and that $H$ is conjugate w.r.t $A$ on $\Lambda$. Then, a LAP holds for $H$ on $\Lambda \setminus \sigma_{\text{pp}}(H)$ w.r.t $A$ and $H$ has no singular continuous spectrum in $\Lambda$.
\end{thm}

For later convenience, we recall that:
\begin{lem}\label{ad-compact} Let ${\cal H}$ be a Hilbert space and $A$ a self-adjoint
operator defined on ${\cal H}$ with domain ${\cal D}(A)$. If $B$ is a compact operator on ${\cal H}$ which belongs to ${\cal C}^{1,1}(A)$, then $\mathrm{ad}_A B$ is also compact.
\end{lem}
The proof of Lemma \ref{ad-compact} corresponds actually to the remark (ii) made in the proof of Theorem 7.2.9 in \cite{abmg}. Due to the inclusions (5.2.10) noted in \cite{abmg}, $\mathrm{ad}_A B$ can be expressed as the norm-limit when $\varepsilon$ tends to 0, of the family of compact operators $(-i\varepsilon^{-1}(e^{iA\varepsilon}B e^{-iA\varepsilon}-B))_{\varepsilon >0}$.

The next result provides a practical criterion to prove the fractional regularity properties mentioned above:
\begin{thm}\label{gsah} Let $Q$ be a self-adjoint operator in ${\cal H}$ bounded from below by a strictly positive constant such that $A^l Q^{-l}$ is continuous for some integer $l\in {\mathbb N}$. Let $0\leq s < l$. Then a bounded symmetric operator $B$ is of class ${\cal C}^{s,1}(A)$ if there exists a function $\chi \in C^{\infty}_0((0,\infty))$ which is positive on some interval $(a,b)$ ($0<a<b<\infty$) such that:
\begin{equation}\label{gsah2}
\int_1^{\infty} \|r^s\chi(Q/r)B\| \frac{dr}{r} < \infty
\end{equation}
\end{thm}
See Theorem 7.5.8 in \cite{abmg} and Theorem 6.1 in \cite{sah} for a proof.

The proof of Theorem \ref{toeplitz} is the result of an interplay between some regularity and conjugacy issues, which are treated in Section 3.2 and 3.3 below.

%*********************
\subsection{Regularity issues}

In the following, we write: ${\cal S}=\bra e_n; n\in {\mathbb N}\ket$.

\begin{lem}\label{bound} Let $F$ be a sesquilinear form defined on ${\cal S}\times {\cal S}$. Assume that
\begin{equation*}
\sum_{(p,q)\in {\mathbb N}^2} |F(e_p,e_q)|^2 < \infty \, .
\end{equation*}
Then, $F$ is continuous on ${\cal S}\times {\cal S}$ for the topology induced by ${\cal H}\times {\cal H}$. It extends continuously to ${\cal H}\times {\cal H}$. If $B$ denotes the (unique) bounded operator associated to that extension, then
\begin{equation*}
\|B\| \leq \left(\sum_{(p,q)\in {\mathbb N}^2} |F(e_p,e_q)|^2 \right)^{1/2} \, .
\end{equation*}
\end{lem}
\pf Let $(\varphi,\psi) \in {\cal S}\times {\cal S}$. We have that:
\begin{equation*}
F(\varphi, \psi)= \sum_{(p,q)\in {\mathbb N}^2} \overline{\bra e_p, \varphi\ket}\bra e_q, \psi\ket F(e_p,e_q) \, .
\end{equation*}
Applying twice Cauchy-Schwarz inequality entails:
\begin{equation*}
|F(\varphi, \psi)|^2 \leq \sum_{(p,q)\in {\mathbb N}^2} |F(e_p,e_q)|^2 \|\varphi\|^2 \|\psi\|^2
\end{equation*}
which proves the first statement. The conclusion is straightforward. \ep

\begin{lem}\label{est} Let $(f,g)\in {\cal A}({\mathbb T})\times {\cal A}({\mathbb T})$ be real-valued  and $(\Phi,\Psi)$ be two complex-valued continuous functions on $(0,\infty)$. Consider the sesquilinear form $F$ defined on ${\cal S}\times {\cal S}$ by: $F(\varphi, \psi)= \bra \varphi, \Phi(X) H_f^* H_g \Psi(X) \psi \ket$. Then, for all $(p,q)\in {\mathbb N}^2$ and all $(\alpha,\beta)\in (1,\infty)^2$ such that $\alpha^{-1}+\beta^{-1}=1,$
\begin{equation*}
|F(e_p,e_q)| \leq |\Phi(p)| |\Psi(q)| \left( \sum_{k\geq p} |\hat{f}_k|^{\alpha} \right)^{1/\alpha}\left( \sum_{k\geq q} |\hat{g}_k|^{\beta} \right)^{1/\beta} \, .
\end{equation*}
\end{lem}
\pf We note that:
\begin{eqnarray*}
F(e_p,e_q)&=& \langle e_p, \Phi(X)H_f^* H_g\Psi(X) e_q \rangle=\overline{\Phi(p)} \Psi(q) \langle e_p, H_f^* H_g e_q \rangle \\
&=& \overline{\Phi(p)} \Psi(q) \sum_{k\leq 0}\langle e_p, H_f^* e_k \ket \bra e_k, H_g e_q \rangle  \, .
\end{eqnarray*}
The conclusion follows from H\"older inequality. \ep

\begin{lem}\label{est1} Let $(f,g)\in C^2({\mathbb T})\times C^2({\mathbb T})$ be two real-valued functions. Consider the sesquilinear forms $F_{\pm}$ and $F$ defined on ${\cal S}\times {\cal S}$ by:
\begin{eqnarray*}
F_+(\varphi, \psi) &=& \bra \varphi, X H_f^* H_g \psi \ket \\
F_-(\varphi, \psi) &=& \bra \varphi, H_g^* H_f X \psi \ket\\
F(\varphi, \psi) &=&  \bra \varphi, X [T_f,T_g]\psi \ket + \bra \varphi, [T_f,T_g]X \psi \ket \,.
\end{eqnarray*}
Then, $F_{\pm}$, $F$ are continuous on ${\cal S}\times {\cal S}$ for the topology induced by ${\cal H}\times {\cal H}$. 
\end{lem}
\pf We prove first the continuity of $F_+$. This is merely a consequence of Lemmata \ref{bound} and \ref{est}, once observed that the Fourier coefficients of the functions $f$ and $g$ satisfy $\sup_n n^2 |\hat{f}_n| <\infty$, $\sup_n n^2 |\hat{g}_n|< \infty$ and by choosing $\alpha >2$ (i.e. $\beta< 2$) in Lemma \ref{est}. Then, we observe that for all $(\varphi,\psi)\in {\cal S}\times {\cal S}$, $F_-(\varphi, \psi)=\overline{F_+(\psi,\varphi)}$. So, the continuity of $F_-$ follows, which implies the continuity of the sesquilinear form $F_s$ defined on ${\cal S}\times {\cal S}$ by: $F_s(\varphi, \psi)=F_+(\varphi, \psi)-F_-(\varphi, \psi)= \bra \varphi, X H_f^* H_g \psi \ket- \bra \varphi, H_g^* H_f X \psi \ket$. We note that the roles of the functions $f$ and $g$ can be exchanged in the above discussion without changing the conclusions. Using formula (\ref{sarason}), we deduce the continuity of $F$.  \ep

The bounded operator associated to the extensions of $F_+$, $F_-$ and $F$ to ${\cal H}\times {\cal H}$ are denoted $X H_f^* H_g$, $H_g^* H_f X$ and $(X [T_f,T_g]+[T_f,T_g]X)$ respectively. We deduce that:
\begin{prop}\label{firstcom} Let $(f,g)\in C^2({\mathbb T})\times C^2({\mathbb T})$ be two real-valued functions. Then, $T_f\in C^1(A_g)$ and
\begin{equation}\label{adgf}
i\mathrm{ad}_{A_g} T_f = \frac{1}{2}( T_g T_{f'} +T_{f'} T_g  )  + \frac{i}{2}([T_g , T_f] X + X [T_g, T_f] )\, .
\end{equation}
\end{prop}
\pf By Lemma \ref{sacom}, the domain ${\cal S}$ is a core for $A_g$. Working with sesquilinear forms on ${\cal S}\times {\cal S}$, we have that:
\begin{equation*}
i(A_g T_f- T_f A_g) = \frac{1}{2}( T_g T_{f'} +T_{f'} T_g  )  + \frac{1}{2}(i[T_g , T_f] X + X i[T_g, T_f] )\, .
\end{equation*}
Due to Lemma \ref{est1}, the RHS is continuous on ${\cal S}\times {\cal S}$ w.r.t. the topology of ${\cal H}\times {\cal H}$. The conclusion follows. \ep

\begin{cor}\label{firstcom1} Let $f \in C^3({\mathbb T})$ be a real-valued function. Then, $T_f$ and $T_{f'}$ belong to $C^1(A_{f'})$. In particular, we have that:
\begin{equation}\label{adf'f}
i\mathrm{ad}_{A_{f'}} T_f = T_{f'}^2 +\frac{i}{2}( [T_{f'} , T_f] X + X [T_{f'}, T_f] ) \, .
\end{equation}
\end{cor}

In order to prove Theorem \ref{toeplitz}, we require a little bit more.
\begin{lem}\label{est2} Let $(f,g)\in C^3({\mathbb T})\times C^2({\mathbb T})$ be two real-valued functions. Consider the sesquilinear forms $G_{\pm}$ defined on ${\cal S}\times {\cal S}$ by:
\begin{eqnarray*}
G_+(\varphi, \psi) &=& \bra \varphi, X H_f^* H_g X \psi \ket \\
G_-(\varphi, \psi) &=& \bra \varphi, X H_g^* H_f X \psi \ket \,.
\end{eqnarray*}
Then, $G_{\pm}$ are continuous on ${\cal S}\times {\cal S}$ w.r.t. the topology induced by ${\cal H}\times {\cal H}$. 
\end{lem}
\pf We observe first that for all $(\varphi,\psi)\in {\cal S}\times {\cal S}$, $G_-(\varphi, \psi)=\overline{G_+(\psi,\varphi)}$. So it is enough to prove the continuity of $G_+$. This is again a consequence of Lemmata \ref{bound} and \ref{est}, once observed that the Fourier coefficients of the functions $f$ and $g$ satisfy $\sup_n n^3 |\hat{f}_n| <\infty$, $\sup_n n^2 |\hat{g}_n|< \infty$ and by choosing $\alpha <2$ (i.e. $\beta> 2$) in Lemma \ref{est}. \ep

The bounded operator associated to the extensions of $G_+$ and $G_-$ to ${\cal H}\times {\cal H}$ are denoted $X H_f^* H_g X$, $X H_g^* H_f X$ respectively. We deduce that:
\begin{prop}\label{compactness} Let $(f,g)\in C^3({\mathbb T})\times C^2({\mathbb T})$ be two real-valued functions. Then, the operator $([T_g, T_f]X + X[T_g, T_f])$ is compact.
\end{prop}
\pf Note that $X$ is invertible ($X\geq 1$) and $X^{-1}$ is compact. Therefore, if $(f,g)\in C^3({\mathbb T})\times C^2({\mathbb T})$, the operators $X H_f^* H_g=(X H_f^* H_g X)X^{-1}$ and $H_f^* H_g X=X^{-1}(X H_f^* H_gX)$ are compact by Lemma \ref{est2}. So are their adjoints, $H_g^* H_fX$ and $X H_g^* H_f$. The conclusion follows since $([T_g,T_f]X+X [T_g,T_f])=X H_f^* H_g-X H_g^* H_f+H_f^* H_g X- H_g^* H_f X$. \ep

\begin{lem}\label{est3} Let $(f,g)\in C^3({\mathbb T})\times C^2({\mathbb T})$ be two real-valued functions. Then, the operators $XH_f^* H_g$, $H_f^* H_g X$, $XH_g^* H_f$ and $H_g^* H_f X$ belong to ${\cal C}^{0,1}(A_g)$. In particular, $([T_g,T_f]X+X [T_g,T_f])$ belongs to ${\cal C}^{0,1}(A_g)$.
\end{lem}
\pf Since the class ${\cal C}^{0,1}(A_g)$ is stable under adjunction $*$, it is enough to prove the result for $H_f^* H_g X$ and $H_g^* H_f X$. With Theorem \ref{gsah} in view, it is enough to show that:
\begin{align*}
& \int_1^{\infty} \| \chi (X/r) H_f^* H_g X\|\, \frac{dr}{r} < \infty \\
& \int_1^{\infty} \| \chi (X/r) H_g^* H_f X\|\, \frac{dr}{r} < \infty \, ,
\end{align*}
where $\chi$ is the characteristic function of some interval $[a,b]\subset (0,\infty)$. Given $0<a<b<\infty$, we deduce from Lemma \ref{est2} that:
\begin{equation*}
\| \chi (X/r) H_f^* H_g X\| \leq \|\chi (X/r) X^{-1}\| \|XH_f^* H_g X\| \leq  \frac{1}{ar}\|XH_f^* H_g X\|
\end{equation*}
which show the finiteness of the first integral. The second case is similar. The last statement follows from (\ref{sarason}). \ep

\begin{prop}\label{TC11} Let $(f,g)\in C^3({\mathbb T})\times C^2({\mathbb T})$ be two real-valued functions. Then, $T_f \in \mathcal{C}^{1,1}(A_g)$. In particular, $T_f \in \mathcal{C}^{1,1}(A_{f'})$.
\end{prop}
\pf It follows from Proposition \ref{firstcom} that $T_{f'}$ and $T_g$ belong to $C^1(A_g)$. In particular, the products $T_{f'}T_g$ and $T_gT_{f'}$ belong also to $C^1(A_g)\subset {\cal C}^{0,1}(A_g)$. By Lemma \ref{est3}, $([T_g, T_f]X + X[T_g, T_f])\in {\cal C}^{0,1}(A_g)$, so we have proven that $\mathrm{ad}_{A_g} T_f$ belongs to ${\cal C}^{0,1}(A_g)$, hence the result. \ep

\noindent{\bf Remark:} A more involved computation shows that if $g\in C^2({\mathbb T})$ and $f\in C^4({\mathbb T})$ are real-valued functions then $T_f \in C^2(A_g)$. Under more restrictive conditions on the symbols $f$ and $g$, the commutators $[T_f, T_g]$ are finite-rank \cite{ding}.

%***********************************************
\subsection{Conjugacy issues}

\noindent{\bf Notations.} Let $f$, $g$ be two real-valued continuous functions defined on ${\mathbb T}$ and $\Lambda \subset {\mathbb R}$ be a Borel set such that $\Lambda \cap$ Ran $f \neq \emptyset$. We define:
\begin{equation*}
c_{\Lambda,f,g}:=\min_{\theta \in \overline{f^{-1}(\Lambda)}}g(\theta) \quad, \quad C_{\Lambda,f,g}:=\max_{\theta \in \overline{f^{-1}(\Lambda)}}g(\theta) \, .
\end{equation*}
If ${\cal I}_{\Lambda}$ denotes the collection of all open sets $\Lambda' \subset {\mathbb R}$ such that $\overline{\Lambda}\subset \Lambda'$, we also write:
\begin{equation*}
c_{\Lambda,f,g}^{\sharp}:=\sup_{\Lambda'\in {\cal I}_{\Lambda}} c_{\Lambda',f,g} \quad, \quad C_{\Lambda,f,g}^{\flat}:=\inf_{\Lambda'\in {\cal I}_{\Lambda}} C_{\Lambda',f,g} \, .
\end{equation*}
\noindent{\bf Remark:} We note that $\overline{f^{-1}(\Lambda)}$ is a compact subset of ${\mathbb T}$. Clearly, $0\leq c_{\Lambda,f,g}^{\sharp} \leq c_{\Lambda,f,g} \leq C_{\Lambda,f,g} \leq C_{\Lambda,f,g}^{\flat}$.

\begin{lem}\label{riez-alg} Let $f\in C^0({\mathbb T})$ be real-valued and $\Phi$ be a complex-valued continuous function vanishing outside the compact set Ran $f$. Then,
\begin{itemize}
\item[(a)] $\Phi(T_f) \simeq P\Phi(L_f) P$
\item[(b)] $P\Phi(L_f)P^{\perp}$ , $P^{\perp}\Phi(L_f)P$ are compact, i.e. $P\Phi(L_f)P^{\perp}\simeq 0 \simeq P^{\perp}\Phi(L_f)P$.
\end{itemize}
\end{lem}
\pf We drop the subscript $f$ and write: $L=L_f$, $T=T_f$. Let us prove by induction on $j$, $j\geq 0$ that $T^j \simeq P L^j P$. This is clear for $j=0$ if we identify the operator $P$ defined on $l^2({\mathbb Z})$ with the identity on $l^2({\mathbb N})$. Assume the induction hypothesis for some $j\geq 0$. Then,
\begin{equation*}
T^{j+1} = (PLP)^{j+1}= (PLP)^j (PLP) \simeq PL^jP(PLP) = PL^j(I-P^{\perp})LP \simeq PL^{j+1}P
\end{equation*}
since $P^{\perp}LP$ is compact by Hartman Theorem. This proves statement (a) for all polynomials $\Phi$. The conclusion follows from Stone-Weierstrass Theorem. Now let us prove by induction on $j$, $j\geq 0$, that $PL^jP^{\perp}\simeq 0$. This is clear for $j=0$. Assuming the induction hypothesis for some $j\geq 0$, we have that:
\begin{equation*}
PL^{j+1}P^{\perp} = (PL^jP)(PLP^{\perp}) + (PL^j P^{\perp})(P^{\perp}LP^{\perp}) \simeq 0
\end{equation*}
by Hartman Theorem ($(PLP^{\perp})$ is compact) and the induction hypothesis ($(PL^j P^{\perp})$ is compact). This proves $P\Phi(L_f)P^{\perp}\simeq 0$ for all polynomials $\Phi$. The conclusion follows again from Stone-Weierstrass Theorem. The proof is complete since $P^{\perp}\Phi(L_f)P=(P\bar{\Phi}(L_f)P^{\perp})^*\simeq 0$. \ep

Another consequence of Stone-Weierstrass Theorem is:
\begin{lem}\label{dif-compact} Let $H_1, H_2$ be two bounded self-adjoint operators defined on $\mathcal{H}$ such that $H_1 \simeq H_2$. For any function $\Phi$ continuous on $\sigma(H_1) \cup \sigma(H_2)$, $\Phi(H_1) \simeq \Phi(H_2)$.
\end{lem}

\begin{prop}\label{Mourre-T} Let $f\in C^3({\mathbb T})$ be a non-constant real-valued symbol. Let $\Lambda \subset$ Ran $f$ be a Borel set. For any real-valued function $\Phi \in C^0({\mathbb R})$ vanishing outside $\Lambda,$
\begin{equation*}
C_{\Lambda,f,|f'|^2}\Phi(T_f)^2 \gtrsim \Phi(T_f)(i\mathrm{ad}_{A_{f'}} T_f)\Phi(T_f) \gtrsim c_{\Lambda,f,|f'|^2}\Phi(T_f)^2
\end{equation*}
If $\overline{\Lambda}\subset$ Ran $f \setminus f(\kappa_f)$, then $c_{\Lambda,f,|f'|^2}>0$.
\end{prop}

\pf Note that $E_{\Lambda}(L_f)=\chi_{f^{-1}(\Lambda)}(L_f)$ for any Borel set $\Lambda \subset {\mathbb R}$. For any $\Lambda \subset$ Ran $f$ and any $\psi \in l^2({\mathbb Z})$, $C_{\Lambda,f,|f'|^2} \| E_{\Lambda}(L_f) \psi  \|^2 \geq \langle  E_{\Lambda}(L_f)\psi,  L_{f'}^{\, 2} E_{\Lambda}(L_f)\psi \rangle \geq c_{\Lambda,f,|f'|^2} \| E_{\Lambda}(L_f) \psi  \|^2$. So, for any real-valued function $\Phi\in C^0({\mathbb R})$ vanishing outside $\Lambda$,
\begin{equation}\label{prec}
C_{\Lambda,f,|f'|^2}\Phi(L_f)^2 \geq \Phi(L_f) L_{f'}^2 \Phi(L_f) \geq c_{\Lambda,f,|f'|^2}\Phi(L_f)^2 \, .
\end{equation}
By Corollary \ref{firstcom1} and Proposition \ref{compactness}, we have that $i\mathrm{ad}_{A_{f'}} T_f \simeq T_{f'}^2$, which implies that:
\begin{eqnarray*}
\Phi(T_f) (i\mathrm{ad}_{A_{f'}} T_f) \Phi(T_f) & \simeq & \Phi(T_f) T_{f'}^2 \Phi(T_f) \\
 &\simeq & P\Phi(L_f) P (PL_{f'}^2P) P\Phi(L_f)P \simeq P\Phi(L_f) L_{f'}^2 \Phi(L_f)P
\end{eqnarray*}
by Lemma \ref{riez-alg}. Since $P\Phi(L_f)^2 P \simeq \Phi(T_f)^2$, the first statement follows from (\ref{prec}). The second statement is a direct consequence of the first one since $f$ is continuous and $\overline{f^{-1}(\Lambda)}\subset f^{-1}(\overline{\Lambda})\subset {\mathbb T}\setminus f^{-1}(f(\kappa_f))\subset {\mathbb T}\setminus \kappa_f$. \ep

The next result shows that the conjugacy property still holds for adequate compact perturbations of $T_f$. We recall that if $V$ is a compact operator, $\sigma_{ess}(T_f+V)=\sigma_{ess}(T_f)=$ Ran $f$.
\begin{cor}\label{Mourre-T3} Let $f\in C^3({\mathbb T})$ be a non-constant real-valued symbol $f$. Let $\Lambda \subset$ Ran $f$ be a Borel set. Let $V$ be a compact symmetric operator defined on $l^2({\mathbb N})$ with $V \in C^1(A_{f'})$ and such that $\mathrm{ad}_{A_{f'}} V$ is compact. Denote $H:=T_f+V$. For any real-valued function $\Phi\in C^0({\mathbb R})$ vanishing outside $\Lambda,$
\begin{equation*}
C_{\Lambda,f,|f'|^2}\Phi(H)^2 \gtrsim \Phi(H) (i\mathrm{ad}_{A_{f'}} H) \Phi(H) \gtrsim c_{\Lambda,f,|f'|^2}\Phi(H)^2
\end{equation*}
If $\overline{\Lambda}\subset$ Ran $f \setminus f(\kappa_f)$, then $c_{\Lambda,f,|f'|^2}>0$.
\end{cor}
\pf Due to Corollary \ref{firstcom1} and the fact that $V \in C^1(A_{f'})$, $H \in C^1(A_{f'})$ and
\begin{equation*}
\Phi(H) (i\mathrm{ad}_{A_{f'}} H) \Phi(H) = \Phi(H) (i\mathrm{ad}_{A_{f'}} T_f) \Phi(H) +  \Phi(H) (i\mathrm{ad}_{A_{f'}} V) \Phi(H) \simeq \Phi(H) (i\mathrm{ad}_{A_{f'}} T_f) \Phi(H) \, .
\end{equation*}
Since $H\simeq T_f$, then $\Phi(H)\simeq \Phi(T_f)$ and $\Phi^2(H)\simeq \Phi^2(T_f)$ by  Lemma \ref{dif-compact}. So, $\Phi(H) (i\mathrm{ad}_{A_{f'}} H) \Phi(H) \simeq \Phi(T_f) (i\mathrm{ad}_{A_{f'}} T_f) \Phi(T_f)$ and the conclusions follow from Proposition \ref{Mourre-T}. \ep

\begin{cor}\label{Mourre-T2} Let $f\in C^3({\mathbb T})$ be a non-constant real-valued symbol $f$. Let $\Lambda \subset$ Ran $f$ be a Borel set. Let $V$ be a compact symmetric operator defined on $l^2({\mathbb N})$ with $V \in C^1(A_{f'})$ and such that $\mathrm{ad}_{A_{f'}} V$ is compact. Denote $H:=T_f+V$. For any open set $\Lambda'$ such that $\overline{\Lambda} \subset \Lambda',$
\begin{equation*}
C_{\Lambda',f,|f'|^2}E_{\Lambda}(H) \gtrsim E_{\Lambda}(H) (i\mathrm{ad}_{A_{f'}} H) E_{\Lambda}(H) \gtrsim c_{\Lambda',f,|f'|^2}E_{\Lambda}(H)\, .
\end{equation*}
If $\overline{\Lambda}\subset$ Ran $f \setminus f(\kappa_f)$, then there exists an open set $\Lambda'$ such that $\overline{\Lambda} \subset \Lambda'$ and $c_{\Lambda',f,|f'|^2}>0$.
\end{cor}
\pf From Corollary \ref{Mourre-T3}, we have that for any open set $\Lambda'$ such that $\overline{\Lambda}\subset \Lambda'$ and any real-valued continuous function $\Phi$ vanishing on ${\mathbb T}\setminus {\Lambda'}$, which takes value $1$ on $\overline{\Lambda}$ (Urysohn Lemma),
\begin{equation*}
C_{\Lambda',f,|f'|^2}\Phi(H)^2 \gtrsim \Phi(H) (i\mathrm{ad}_{A_{f'}} H) \Phi(H) \gtrsim c_{\Lambda',f,|f'|^2}\Phi(H)^2
\end{equation*}
The first statement follows after multiplying the previous inequalities on both sides by $E_{\Lambda}(H)$. Given $\Lambda$ such $\overline{\Lambda}\subset$ Ran $f \setminus f(\kappa_f)$, we can pick an open set $\Lambda'$ such that $\overline{\Lambda} \subset \Lambda'$ and $\overline{\Lambda'}\subset$ Ran $f \setminus f(\kappa_f)$ (Ran $f \setminus f(\kappa_f)$ is an open subset of Ran $f$). Due to Corollary \ref{Mourre-T3}, $0< c_{\Lambda',f,|f'|^2}$. \ep

%****************************************
\subsection{Proof of Theorem \ref{toeplitz}}

By hypothesis $V\in \mathcal{C}^{1,1}(A_{f'})$ is compact, thus $\mathrm{ad}_{A_{f'}} V$ is also compact (see Lemma \ref{ad-compact}). By Proposition \ref{TC11}, $T_f\in \mathcal{C}^{1,1}(A_{f'})$. So, $H=T_f+V \in \mathcal{C}^{1,1}(A_{f'})$. Due to Corollary \ref{Mourre-T2}, $H$ is conjugate w.r.t $A_{f'}$ on any Borel set and any open set $\Lambda$ such that $\overline{\Lambda}\subset$ Ran $f \setminus f(\kappa_f)$. Statement (a) follows from Proposition \ref{virial-2} while statement (b) derives from Theorem \ref{nosc2}.

%****************************************
\subsection{Proof of Theorem \ref{schro}}

The main point consists in proving that the operator $V$ belongs to ${\cal C}^{1,1}(A_{f'})$. Let $S$ denote the unilateral shift: $Se_n =e_{n+1}$, $n\in {\mathbb N}$. Note that $S^* e_1=0$ and $S^* e_n = e_{n-1}$ otherwise. This allows us to rewrite: $T_f=S+S^*$ and $T_{f'}=i(S-S^*)$. Note that $S$ and $S^*$ belong to $C^{\infty}(X)$. In particular, $\mathrm{ad}_X S=S$ and $\mathrm{ad}_X S^*=-S^*$. We have that:
\begin{equation*}
A_{f'}= \frac{i}{2}\left((S-S^*) X + X (S-S^*)\right) \, .
\end{equation*}
We use the following local notations. To any bounded sequence $\gamma:=(\gamma_k)_{k\in {\mathbb N}}$ in ${\mathbb C}^{\mathbb N}$, we associate the bounded linear operator $D_{\gamma}$ defined by its action on the canonical orthonormal basis of $l^2({\mathbb N})$: $D_{\gamma} e_n=\gamma_n e_n$, $n\in {\mathbb N}$. We recall that $\|D_{\gamma}\|=\sup_{n}|\gamma_n|=q_0(\gamma)$. Note that if $\gamma$ and $\beta$ are two bounded sequences in ${\mathbb C}^{\mathbb N}$ and $c\in {\mathbb C}$, then $D_{\gamma+\beta}=D_{\gamma}+D_{\beta}$, $D_{\gamma\cdot \beta}=D_{\gamma}D_{\beta}$, $D_{c\gamma}=cD_{\gamma}$, $D_{\gamma}^*=D_{\bar{\gamma}}$ and $[D_{\gamma},D_{\beta}]=0$. In addition, if $\xi\cdot\gamma$ is bounded, then $D_{\xi\cdot \gamma}=XD_{\gamma}=D_{\gamma}X$. For any bounded sequence $\gamma$, $D_{\gamma} \in C^{\infty}(X)$ and $\mathrm{ad}_X D_{\gamma}=0$.

\begin{lem}\label{criteria1} If $q_1(\gamma)<\infty$, then $D_{\gamma}\in C^1(A_{f'})$.
\end{lem}
\pf As a sesquilinear form on ${\cal D}_X\times {\cal D}_X$, we have that:
\begin{equation*}
i(A_{f'}D_{\gamma}-D_{\gamma}A_{f'})= -\frac{1}{2} \left([S,D_{\gamma}]X-[S^*,D_{\gamma}]X +X [S,D_{\gamma}]-X[S^*,D_{\gamma}]\right)
\end{equation*}
where $[S,D_{\gamma}]= D_{\Delta \gamma}S$ and $[S^*,D_{\gamma}]=-S^* D_{\Delta \gamma}$. If $q_1(\gamma)<\infty$, the identity extends by continuity to ${\cal H}\times {\cal H}$. Since ${\cal D}_X$ is a core for $A_{f'}$, $D_{\gamma}\in C^1(A_{f'})$ and
\begin{equation}\label{expl}
i\mathrm{ad}_{A_{f'}} D_{\gamma}= -\frac{1}{2} \left(D_{2\xi \Delta\gamma-\Delta\gamma}S + S^* D_{2\xi \Delta\gamma-\Delta\gamma}\right) \, .
\end{equation} \ep

\begin{lem}\label{interm} $S\in C^1(A_{f'})$.
\end{lem}
\pf Note first that $S^*S=I$, $SS^*=I-|e_1\ket \bra e_1|$, so $[S^*,S]=|e_1\ket \bra e_1|$. As a consequence, the operators $[S^*,S]X$ and $X[S^*,S]$ (defined via their associated sesquilinear forms on ${\cal D}_X\times {\cal D}_X$) are bounded and $X[S^*,S]=|e_1\ket \bra e_1|=[S^*,S]X$. Using sesquilinear form on ${\cal D}_X\times {\cal D}_X$, we have that:
\begin{equation*}
i(A_{f'}S-SA_{f'})= -\frac{1}{2} \left(S [X,S]-S^*[X,S]-[S^*,S]X+[X,S]S-[X,S]S^*-X[S^*,S]\right)
\end{equation*}
Since ${\cal D}_X$ is a core for $A_{f'}$, the result follows and
\begin{equation*}
i\mathrm{ad}_{A_f'} S = -S^2 +1 +\frac{1}{2} |e_1\ket\bra e_1 |\, .
\end{equation*}
\ep

\begin{lem}\label{criteria2} If $q_2(\gamma)<\infty$, then $D_{\gamma}\in C^2(A_{f'})$.
\end{lem}
\pf Since $q_1(\gamma)\leq q_2(\gamma) < \infty$, $D_{\gamma}\in C^1(A_{f'})$ by Lemma \ref{criteria1}. It remains to prove that $\mathrm{ad}_{A_{f'}} D_{\gamma} \in C^1(A_{f'})$. By Lemma \ref{interm}, $S$ (and $S^*$) belongs to $C^1(A_{f'})$.  Since $q_2(\gamma)< \infty$, then $q_1(2\xi \Delta\gamma-\Delta\gamma)<\infty$ and so $D_{2\xi \Delta\gamma-\Delta\gamma}\in C^1(A_{f'})$ by Lemma \ref{criteria1}. In view of (\ref{expl}), we deduce that $i\mathrm{ad}_{A_{f'}} D_{\gamma} \in C^1(A_{f'})$, hence the result. \ep

\begin{lem}\label{d1approx} Let $\gamma$ be a bounded sequence of real numbers.
\begin{itemize}
\item[(a)] If
\begin{equation}\label{h1}
\int_1^{\infty} \sup_{a_1 r \leq n \leq b_1 \, r} |\gamma_n|\, dr < \infty
\end{equation}
for some $0<a_1 < b_1 < \infty$, then $D_{\gamma} \in {\cal C}^{1,1}(A_{f'})$.
\item[(b)] If $q_1(\gamma)<\infty$ and
\begin{equation}\label{h2}
\int_1^{\infty} \sup_{a_2 r \leq n \leq b_2 r} |\gamma_{n+1} -\gamma_n|\, dr < \infty
\end{equation}
for some $0<a_2<b_2<\infty$, then $D_{\gamma} \in C^1(A_{f'})$ and $\mathrm{ad}_{A_{f'}} D_{\gamma} \in {\cal C}^{0,1}(A_{f'})$. In particular, $D_{\gamma} \in {\cal C}^{1,1}(A_{f'})$.
\end{itemize}
\end{lem}
\noindent{\bf Proof:} The proof follows from Lemma \ref{relbound} and Theorem \ref{gsah} where the roles of $A$ and $Q$ are endorsed by $A_{f'}$ and $X$ respectively. We start with statement (a). Let $\chi$ be a smoothed characteristic function supported on the interval $(a_1,b_1)$:
\begin{equation*}
\int_1^{\infty} \|\chi(X/r)D_{\gamma}\| \, dr \leq  \int_1^{\infty} \sup_{a_1 r \leq n \leq b_1 \, r} |\gamma_n|\, dr < \infty \, ,
\end{equation*}
and the conclusion follows from Theorem \ref{gsah}. We go on with statement (b). Since $q_1(\gamma)<\infty$, $D_{\gamma} \in C^1(A_{f'})$ by Lemma \ref{criteria1}. Let $\chi$ be a smoothed characteristic function supported on the interval $(a_2,b_2)$. Since $\|S\|=1$, we have that
\begin{equation*}
\int_1^{\infty} \|\chi(X/r) D_{2\xi \Delta\gamma-\Delta\gamma}S \| \, \frac{dr}{r} \leq \int_1^{\infty} \|\chi(X/r) D_{2\xi \Delta\gamma-\Delta\gamma}\| \, \frac{dr}{r} \leq C \int_1^{\infty} \sup_{a_2 r \leq n \leq b_2 \, r} |\gamma_{n+1}-\gamma_n |\, dr < \infty \, ,
\end{equation*}
for some $C>0$. By Theorem \ref{gsah}, we get that $D_{2\xi \Delta\gamma-\Delta\gamma}S \in  {\cal C}^{0,1}(A_{f'})$. Since ${\cal C}^{0,1}(A_{f'})$ is stable under adjunction $*$, we deduce that $S^* D_{2\xi \Delta\gamma-\Delta\gamma}\in {\cal C}^{0,1}(A_{f'})$. In view of (\ref{expl}), this shows that $\mathrm{ad}_{A_{f'}} D_{\gamma}$ belong to ${\cal C}^{0,1}(A_{f'})$. \ep

Now, consider the operator $V$ defined in Theorem \ref{schro}. In view of Lemmata \ref{criteria2} and \ref{d1approx}, it is expressed as a sum of operators which belong to ${\cal C}^{1,1}(A)$ (since $C^2(A)\subset {\cal C}^{1,1}(A)$). Therefore, $V\in {\cal C}^{1,1}(A)$ and Theorem \ref{schro} is a consequence of Theorem \ref{toeplitz}.

%****************************************
\subsection{Proof of Theorem \ref{rank1}}

The operators $(V_{\beta})_{\beta\in {\mathbb R}^N}$ are finite rank, so $\sigma_{ess}(H_{\beta})=$ Ran $f$ for all $\beta\in {\mathbb R}^N$. The main point consists in proving that these operators belong to ${\cal C}^{1,1}(A_{f'})$.

\begin{lem}\label{rank11} Let $g\in C^2({\mathbb T})$ be a real-valued function. Let $\varphi$, $\psi$ be two unitary vectors of $l^2({\mathbb{N}})$ and assume that one of them satisfies the hypothesis {\bf H} for some $0<a<b<\infty$. Then, the projectors $|\psi\ket \bra \varphi|$ and $|\varphi\ket \bra \psi|$ belong to $\mathcal{C}^{1,1}(A_g)$.
\end{lem}
\pf In order to fix ideas, assume that $\psi$ satisfies the hypothesis ${\bf H}$. We note that the class ${\cal C}^{1,1}(A_g)$ is stable under adjunction $*$. So it is enough to prove that $|\psi\ket \bra \varphi|$ belongs to $\mathcal{C}^{1,1}(A_g)$. Let $\chi$ be the characteristic function of the interval $[a,b]\subset (0,\infty)$. Then,
$$
\int_{1}^{\infty} \| \chi (X/r) |\psi\ket \bra \varphi| \| \,  dr = \int_{1}^{\infty}
( \sum_{n\in {\mathbb N}\cap [ar, br]} | \psi_n |^2 \, )^{1/2}\, dr \, < \infty
$$
In view of Lemma \ref{relbound} and Theorem \ref{gsah}, one has that $|\psi\ket \bra \varphi| \in \mathcal{C}^{1,1}(A_g)$, hence the conclusion. \ep

\begin{lem}\label{rank12} Let $A$ be a self-adjoint operator defined on some Hilbert space ${\cal H}$, with domain ${\cal D}(A)$. Assume that $\varphi$, $\psi$ are two unitary vectors of ${\cal H}$ which belong to ${\cal D}(A^k)$ for some $k\in {\mathbb N}$. Then, $|\psi\ket \bra \varphi| \in C^k(A)$.
\end{lem}
\pf It is enough to prove by induction on $j$, $j\in \{1,\ldots,k\}$ that $|\psi\ket \bra \varphi| \in C^j(A)$ and
\begin{equation*}
\text{ad}_{A}^j(|\psi\ket \bra \varphi|) = \sum_{p=0}^j (-1)^{j-p}{j \choose p} |A^p\psi><A^{j-p}\varphi | \, .
\end{equation*}
The details are omitted. \ep

Now, consider the operators $(V_{\beta})$ defined in Theorem \ref{rank1}. In view of Lemmata \ref{rank11} and \ref{rank12}, they are expressed as a sum of operators which belong to ${\cal C}^{1,1}(A)$ (since $C^2(A)\subset {\cal C}^{1,1}(A)$). Therefore, $V_{\beta}\in {\cal C}^{1,1}(A)$ for all $\beta\in {\mathbb R}^N$ and Theorem \ref{rank1} is a consequence of Theorem \ref{toeplitz}.

%****************************************
\subsection{Proof of Theorem \ref{product}}

Since $f$ and $g$ belong to $C^3({\mathbb {T}})$, $h=fg \in C^3({\mathbb {T}})$ and the operators $T_f$, $T_g$ and $T_h=T_{fg}$ belong to ${\mathcal C}^{1,1}(A_{h'})$. We also deduce that $T_fT_g$, $T_gT_f$ and $\Re(T_fT_g)$ belong to ${\mathcal C}^{1,1}(A_{h'})$ (see e.g. Proposition 5.2.3 in \cite{abmg}). So, the compact operators $H_f^* H_g = T_h-T_fT_g$ and $H_g^* H_f = T_h-T_gT_f$ also belong to ${\mathcal C}^{1,1}(A_{h'})$. Since
\begin{equation*}
H= T_h +\frac{1}{2}(H_g^* H_f-H_f^* H_g)
\end{equation*}
Theorem \ref{product} is a consequence of Theorem \ref{toeplitz}.

%*****************************
\section{Propagation estimates}

The commutator formalism also allows to derive some propagation estimates. In this section, $f$ denotes a real-valued symbol in $C^3({\mathbb T})$. We characterize first the propagation properties of the Toeplitz operators. Through this section we define for all $\psi \in {\cal D}_X$, $\|\psi\|_X= \sqrt{\|\psi\|^2+\| X \psi \|^2}$.
\begin{prop}\label{propa2} Consider a non-constant real-valued symbol $f \in C^3({\mathbb T})$. Then, for all $\varphi \in \mathcal{D}_X,$
\begin{eqnarray*}
\left( \liminf_{t \to \pm \infty}\frac{1}{t} \int_0^{t} \| T_{f'} e^{iT_f s}\varphi \|^2 \, ds \right)^{1/2} & \leq & \liminf_{t \to \pm \infty}\frac{1}{t} \| e^{iT_f t}\varphi \|_X \\
\limsup_{t \to \pm \infty}\frac{1}{t} \| e^{iT_f t}\varphi \|_X &\leq & \left( \limsup_{t \to \pm \infty}\frac{1}{t} \int_0^{t} \| T_{f'} e^{iT_f s}\varphi \|^2 \, ds \right)^{1/2}
\end{eqnarray*}
\end{prop}
See Section 4.3 for the proof. These properties are preserved in a weaker form under perturbations:
\begin{prop}\label{propa4} Consider a non-constant real-valued symbol $f \in C^3({\mathbb T})$. Let $V\in C^1(A_{f'})\cap C^1(X) \cap C^1(X^2)$ be a bounded symmetric operator on $l^2({\mathbb N})$. Let $H=T_f+V$. Then, for all $\varphi \in {\cal D}_X$
$$
\limsup_{t \to \pm \infty} \frac{1}{|t|} \| e^{iHt}\varphi \|_X \leq \sqrt{\|\mathrm{ad}_{A_{f'}}H \|} \|\varphi \| \, .
$$
\end{prop}
See Section 4.4 for the proof.

\begin{prop}\label{propa3} Consider a non-constant real-valued symbol $f \in C^3({\mathbb T})$. Let $V\in C^1(A_{f'})\cap C^1(X) \cap C^1(X^2)$ be a compact symmetric operator such that $\mathrm{ad}_{A_{f'}}V$ is also compact. Let $H=T_f+V$ and $\Lambda \subset$ Ran $f$ be a Borel set. Then,
\begin{itemize}
\item For all $\varphi \in {\cal H}$ such that $E_{\Lambda}(H)\varphi \in {\cal H}_c(H)\cap {\cal D}_X$
$$
\sqrt{c_{\Lambda,f,|f'|^2}^{\sharp}} \|E_{\Lambda}(H)\varphi \| \leq \liminf_{t \to \pm \infty} \frac{1}{|t|} \| e^{iHt}E_{\Lambda}(H)\varphi \|_X \leq \limsup_{t \to \pm \infty} \frac{1}{|t|} \| e^{iHt}E_{\Lambda}(H)\varphi \|_X \leq \sqrt{C_{\Lambda,f,|f'|^2}^{\flat}} \|E_{\Lambda}(H)\varphi \| \, .
$$
\item For all $\varphi \in {\cal H}_c(H)\cap {\cal D}_X$ and all real-valued $\Phi\in C^{\infty}_0({\mathbb R})$ vanishing outside $\Lambda,$ 
$$
\sqrt{c_{\Lambda,f,|f'|^2}} \|\Phi(H)\varphi \| \leq \liminf_{t \to \pm \infty} \frac{1}{|t|} \| e^{iHt}\Phi(H)\varphi \|_X \leq \limsup_{t \to \pm \infty} \frac{1}{|t|} \| e^{iHt}\Phi(H)\varphi \|_X \leq \sqrt{C_{\Lambda,f,|f'|^2}} \|\Phi(H)\varphi \| \, .
$$
\end{itemize}
If $\overline{\Lambda} \subset$ Ran $f\setminus f(\kappa_f)$, then $c_{\Lambda,f,|f'|^2}\geq c_{\Lambda,f,|f'|^2}^{\sharp} >0$.
\end{prop}
See Section 4.5 for the proof.

\noindent{\bf Remark:} The conclusions of Propositions \ref{propa4} and \ref{propa3} still hold if instead of the condition $V\in C^1(A_{f'})\cap C^1(X) \cap C^1(X^2)$, we only require that $V\in C^1(A_{f'})\cap C^1(X)$ and that for all $\varphi \in {\cal D}_X$
\begin{equation*}
\lim_{t\rightarrow \pm \infty} \frac{1}{t} \Re \langle i(\mathrm{ad}_X V) e^{iH t} \varphi, X e^{iH t} \varphi \rangle =0 \, .
\end{equation*}

\noindent{\bf Remark:} The construction of a non trivial vector $\varphi$ satisfying the condition $E_{\Lambda}(H)\varphi \in {\cal H}_c(H)\cap {\cal D}_X$ under the hypotheses of Proposition \ref{propa3} can be performed as follows. Let $\Lambda$ be an open interval such that $\overline{\Lambda} \subset$ Ran $f\setminus f(\kappa_f)$. Since $V\in C^1(A_{f'})$ is symmetric, we deduce from Corollary \ref{Mourre-T2} and Proposition \ref{virial-2} that $\sigma_{pp}(H)\cap \Lambda$ is finite. Moreover $V\in C^1(X)$, so $H$ and $\Phi(H)$ belong to $C^1(X)$ for any $\Phi \in C^{\infty}_0(\Lambda)$. By considering $\Phi \in C^{\infty}_0(\Lambda \setminus \sigma_{pp}(H))$ and defining $\varphi=\Phi(H)\psi$ where $\psi \in {\cal D}_X$, we have that $E_{\Lambda}(H)\varphi=\varphi \in {\cal D}_X \cap {\cal H}_c(H)$.

Sections 4.1 and 4.2 articulate Mourre Theory with the proofs of Propositions \ref{propa2}, \ref{propa4} and \ref{propa3}.

%******************************************
\subsection{Preliminaries}

Let ${\cal H}$ be a Hilbert space and $A$ be a self-adjoint operator defined on ${\cal H}$ with domain ${\cal D}(A)$. We recall that if $H \in C^1(A)$, then $e^{iHt} \in C^1(A)$ and $e^{iHt} \mathcal{D}(A) \subset \mathcal{D}(A)$ for any $t \in {\mathbb R}$. So given $t \in {\mathbb R}$, let us define on $\mathcal{D}(A) \times \mathcal{D}(A)$ the sequilinear form:
$$
F_t(\varphi, \psi) =  \langle e^{iHt} \varphi, A e^{iHt}\psi \rangle  -  \langle  \varphi, A\psi \rangle = \int_0^t \langle  e^{iHs} \varphi, i [A,H]e^{iHs}\psi \rangle \, ds
$$
$F_t$ is continuous on $\mathcal{D}(A) \times \mathcal{D}(A) $ w.r.t. the topology induced by $\mathcal{H}\times \mathcal{H}$: actually we have that:
$
|F_t(\varphi, \psi )| \leq \|\mathrm{ad}_A H\| |t| \, \|  \varphi\| \,  \| \psi \|\, .
$
If $F_t^o$ denotes the continuous extension of $F_t$ to $\mathcal{H}\times \mathcal{H}$, then for all $(\varphi, \psi)\in {\cal H}\times {\cal H},$
\begin{equation}\label{intform1}
F_t^o(\varphi, \psi )  =  \int_0^t \langle  e^{iHs} \varphi, i(\mathrm{ad}_A H) e^{iHs}\psi \rangle \, ds
\end{equation}
and for all $(\varphi, \psi)\in {\cal H}\times {\cal D}(A)$, $F_t^o(\varphi, \psi )  =  \langle e^{iHt} \varphi, A e^{iHt}\psi \rangle  -  \langle  \varphi, A\psi \rangle$. Thus, for all $\varphi \in {\cal D}(A)$ and all $t\in {\mathbb R},$
\begin{equation}\label{intform2}
e^{-iHt}A e^{iHt}\varphi - A\varphi= \int_0^t e^{-iHs} i(\mathrm{ad}_A H) e^{iHs}\varphi \, ds \, .
\end{equation}
It follows that:
\begin{lem}\label{brut} Let $H\in C^1(A)$. Then, for all $(\varphi,\psi) \in {\cal H}\times {\cal D}(A),$
\begin{equation*}
\limsup_{t\rightarrow \pm\infty} \frac{1}{|t|}  | \langle e^{iHt}\varphi,  Ae^{iHt}\psi \rangle | \leq \|\mathrm{ad}_A H \| \, \|\varphi  \| \,\|\psi  \| \, .
\end{equation*}
\end{lem}

RAGE Theorem \cite{cfks} tells us that for all $\varphi \in \mathcal{H}_c(H)$ and any compact operator $K$ defined on ${\cal H},$
\begin{equation*}
\lim_{t \to \pm \infty} \frac{1}{t}  \int_0^{t} \|Ke^{iHs} \varphi  \| ds= 0\, .
\end{equation*}

We deduce the two following lemmata:
\begin{lem}\label{obs} Let $H\in C^1(A)$ and assume that $i(\mathrm{ad}_A H) \gtreqless B + K$ where $B$ is bounded and $K$ is compact. Then, for all $\varphi \in \mathcal{D}(A)\cap \mathcal{H}_c(H)$
\begin{eqnarray*}
\liminf_{t \to \pm \infty} \, \frac{1}{t}  \langle e^{iHt} \varphi, Ae^{iHt}\varphi \rangle
&\gtreqless & \liminf_{t \to \pm \infty} \frac{1}{t} \int_0^{t}  \langle e^{iHs} \varphi,  Be^{iHs}\varphi \rangle \, ds\\
\limsup_{t \to \pm \infty} \, \frac{1}{t} \langle e^{iHt} \varphi, Ae^{iHt}\varphi \rangle &\gtreqless & \limsup_{t \to \pm \infty} \frac{1}{t} \int_0^{t}  \langle e^{iHs} \varphi,  Be^{iHs}\varphi \rangle \, ds \, .
\end{eqnarray*}
\end{lem}

\begin{lem}\label{mourre-propa} Let $H\in C^1(A)$ and $\Phi \in L^{\infty}({\mathbb R})$ be real-valued.
\begin{itemize}
\item If $\Phi (H) (i\mathrm{ad}_A H) \Phi (H) \geq c\Phi (H)^2 + K$ for some $c\in {\mathbb R}$ and $K$ compact then for any $\varphi \in {\cal H}$ such that $\Phi (H)\varphi \in \mathcal{D}(A) \cap \mathcal{H}_c(H)$, one has that
\begin{equation*}
\liminf_{t\rightarrow \pm\infty} \frac{1}{t} \langle e^{iHt}\Phi (H) \varphi,  Ae^{iHt}\Phi (H) \varphi \rangle \geq c\| \Phi (H)\varphi  \|^2
\end{equation*}
\item If $\Phi (H) (i\mathrm{ad}_A H) \Phi (H) \leq C \Phi (H) + K$ for some $C\in {\mathbb R}$ and $K$ compact then for any $\varphi \in {\cal H}$ such that $\Phi (H)\varphi \in \mathcal{D}(A) \cap \mathcal{H}_c(H)$, one has that
\begin{equation*}
\limsup_{t\rightarrow \pm\infty} \frac{1}{t} \langle e^{iHt}\Phi (H) \varphi,  Ae^{iHt}\Phi (H) \varphi \rangle \leq C\| \Phi (H)\varphi  \|^2
\end{equation*}
\end{itemize}
\end{lem}

Note that the conclusions of Lemma \ref{obs} can be strengthened in the following cases:
\begin{lem}\label{obs2} Let $H\in C^1(A)$. Assume that $i(\mathrm{ad}_A H)= B + K$ where $B$ is bounded, $K$ is compact and $[B,H]=0$.
\begin{itemize}
\item[(a)] If $K=0$, then for all $\varphi \in \mathcal{D}(A),$ \begin{equation*}
\lim_{t\rightarrow \pm\infty} \frac{1}{t}  e^{-iHt}Ae^{iHt}\varphi = i(\mathrm{ad}_A H)\varphi =B\varphi \, .
\end{equation*}
\item[(b)] If $\sigma_{pp}(H)$ is finite, then for all $\varphi \in \mathcal{D}(A),$
\begin{equation*}
\lim_{t\rightarrow \pm\infty} \frac{1}{t} e^{-iHt}Ae^{iHt}\varphi = E_c(H) BE_c(H)\varphi\, .
\end{equation*}
\end{itemize}
\end{lem}
\noindent{\bf Proof:} Case (a) follows from (\ref{intform2}). Now, consider Case (b). Assume that $\sigma_{pp}(H)\neq \emptyset$, denote by $(\lambda_j)_{j=1}^N$, $N\in {\mathbb N}$ the set of eigenvalues of $H$ and by $E_{pp}(H)$ (resp. $E_c(H)$) the orthogonal projection on ${\cal H}_{pp}(H)$ (resp. ${\cal H}_c(H)$). For all $\varphi \in {\cal H}$ and all $t\neq 0$, we have that:
\begin{equation*}
\frac{1}{t}\int_0^t e^{-iHs} K e^{iHs} E_{pp}(H) \varphi \, ds = \sum_{j=1}^N \frac{1}{t}(\int_0^t e^{-i(H-\lambda_j)s} \, ds ) K E_{\{\lambda_j\}}(H) \varphi
\end{equation*}
which tends to $\sum_{j=1}^N E_{\{\lambda_j\}}(H) K E_{\{\lambda_j\}}(H) \varphi$ when $t$ tends to $\pm \infty$ (see e.g. Theorem 1.3 in \cite{kato} Chapter X). By RAGE Theorem, we have also that for all $\varphi \in {\cal H},$
\begin{equation*}
\lim_{t \rightarrow \pm \infty} \frac{1}{t} \int_0^t e^{-iHs} K e^{iHs} E_c(H) \varphi \, ds =0
\end{equation*}
In other words, using the fact that $[B,H]=0$, we get that for all $\varphi \in {\cal H},$
\begin{equation*}
B\varphi + \lim_{t \rightarrow \pm \infty} \frac{1}{t} \int_0^t e^{-iHs} K e^{iHs}\varphi \, ds = \sum_{j=1}^N E_{\{\lambda_j\}}(H)(B +K) E_{\{\lambda_j\}}(H)\varphi + E_c(H)B E_c(H)\varphi
\end{equation*}
But, we have that $E_{\{\lambda_j\}}(H)i[A,H] E_{\{\lambda_j\}}(H)=E_{\{\lambda_j\}}(H)(B +K) E_{\{\lambda_j\}}(H)=0$ for all $j\in \{1,\ldots,N\}$, due to the Virial Theorem (see Proposition 7.2.10 in \cite{abmg}). The conclusion follows from (\ref{intform2}). If $\sigma_{pp}(H)= \emptyset$, we proceed as before, noting that $E_{pp}(H)=0$ (and $E_c(H)=I$). \ep

%******************************************
\subsection{Preliminaries. continued}

We come back to the context and notations of Section 2. Recall that if $H \in C^1(X)$, then $e^{iHt}$ and $He^{iHt} \in C^1(X)$ for any $t \in \mathbb{R}$. So, $e^{iHt} \mathcal{D}_X \subset \mathcal{D}_X$ and $He^{iHt}\mathcal{D}_X \subset \mathcal{D}_X$ for all $t \in {\mathbb R}$.

\begin{lem}\label{propa1} Let $H \in C^1 (X)$. Then for all $t \in \mathbb{R}$ and all $\varphi \in {\cal D}_X,$
$$
\| e^{iHt}\varphi \|_X^2 - \| \varphi \|_X^2 = 2 \int_0^t \Re (\langle i(\mathrm{ad}_X H) e^{iHs} \varphi, X e^{iHs} \varphi \rangle \, ds\, .
$$
As a consequence, if $\limsup_{t\rightarrow \pm \infty} |t|^{-1}|\Re (\langle i(\mathrm{ad}_X H) e^{iHt} \varphi, X e^{iHt} \varphi \rangle| <\infty$ then,
\begin{equation*}
\limsup_{t\rightarrow \pm \infty} \frac{1}{|t|} \| e^{iHt}\varphi \|_X \leq \left(\limsup_{t\rightarrow \pm \infty} \frac{1}{|t|}|\Re (\langle i(\mathrm{ad}_X H) e^{iHt} \varphi, X e^{iHt} \varphi \rangle| \right)^{1/2} < \infty
\end{equation*}
If $0 \leq \liminf_{t\rightarrow \pm \infty} t^{-1}\Re (\langle i(\mathrm{ad}_X H) e^{iHt} \varphi, X e^{iHt} \varphi \rangle \leq \limsup_{t\rightarrow \pm \infty} t^{-1} \Re (\langle i(\mathrm{ad}_X H) e^{iHt} \varphi, X e^{iHt} \varphi \rangle <\infty$, then
\begin{eqnarray*}
\left(\liminf_{t\rightarrow \pm \infty} \frac{1}{t}\Re (\langle i(\mathrm{ad}_X H) e^{iHt} \varphi, X e^{iHt} \varphi \rangle \right)^{1/2} &\leq & \liminf_{t\rightarrow \pm \infty} \frac{1}{t} \| e^{iHt}\varphi \|_X \\
\text{and}\quad \limsup_{t\rightarrow \pm \infty} \frac{1}{t} \| e^{iHt}\varphi \|_X &\leq & \left(\limsup_{t\rightarrow \pm \infty} \frac{1}{t}\Re (\langle i(\mathrm{ad}_X H) e^{iHt} \varphi, X e^{iHt} \varphi \rangle \right)^{1/2} < \infty \, .
\end{eqnarray*}
\end{lem}
\pf We have that for all $t\in {\mathbb R}$ and all $\varphi \in {\cal D}_X$: $\| e^{iHt}\varphi \|_X^2 -\| \varphi \|_X^2 =\| Xe^{iHt}\varphi \|^2- \| X\varphi \|^2$ and so
\begin{equation*}
\| e^{iHt}\varphi \|_X^2 -\| \varphi \|_X^2 = \int_{0}^t \partial_s  \langle X e^{iHs}\varphi, Xe^{iHs}  \varphi \rangle \, ds = 2 \int_0^t  \Re \langle i[X, H] e^{iHs} \varphi, Xe^{iHs} \varphi  \rangle \, ds \, .
\end{equation*}
\ep

\noindent{\bf Remark:} Let $H \in C^1 (X)$ and assume that $[H, \mathrm{ad}_X H]=0$. Then for all $t \in \mathbb{R}$ and all $\varphi \in {\cal D}_X,$
\begin{equation*}
\| e^{iHt}\varphi \|_X^2 - \| \varphi \|_X^2 = t^2 \| (\mathrm{ad}_X H)\varphi \|^2+2t \Re \bra i(\mathrm{ad}_X H)\varphi,X\varphi \ket \, .
\end{equation*}

%******************************************
\subsection{Proof of Proposition \ref{propa2}}

First, we note that $\sigma_{pp}(T)=\emptyset$ \cite{hw}, \cite{ros}. We also have that $T_f\in C^1(A_{f'})$ and $i\mathrm{ad}_{A_{f'}}T_f \simeq T_{f'}^2$ by Corollary \ref{firstcom1} and Proposition \ref{compactness}. We deduce from Lemma \ref{obs} that for all $\varphi \in \mathcal{D}(A_{f'}),$
\begin{eqnarray*}
0\leq \liminf_{t\rightarrow \pm\infty} \frac{1}{t} \int_0^{t} \|T_{f'} e^{iT_f s} \varphi \|^2\, ds &=& \liminf_{t\rightarrow \pm\infty} \frac{1}{t} \langle e^{iT_f t}\varphi,  A_{f'}e^{iT_f t}\varphi \rangle \\
\limsup_{t\rightarrow \pm\infty} \frac{1}{t} \langle e^{iT_f t}\varphi,  A_{f'}e^{iT_f t}\varphi \rangle &=& \limsup_{t\rightarrow \pm\infty} \frac{1}{t} \int_0^{t} \|T_{f'} e^{iT_f s} \varphi \|^2 \, ds < \infty \, .
\end{eqnarray*}
The conclusion follows from Lemma \ref{propa1} by noting that: $i\mathrm{ad}_X T_f =T_{f'}$ and so that $\Re (\langle i(\mathrm{ad}_X T_f) e^{iT_f t} \varphi, X e^{iT_f t} \varphi \rangle= \langle e^{iT_f t} \varphi, A_{f'} e^{iT_f t} \varphi \rangle$ for all $t\in {\mathbb R}$ and all $\varphi \in {\cal D}_X \subset {\cal D}(A_{f'})$.

%******************************************
\subsection{Proof of Proposition \ref{propa4}}

Due to Lemma \ref{ckx} and Corollary \ref{firstcom1}, we have that $H \in C^1(A_{f'})\cap C^1(X)$. By Lemma \ref{brut}, for all $(\varphi,\psi) \in {\cal H}\times {\cal D}(A),$
\begin{equation*}
\limsup_{t\rightarrow \pm\infty} \frac{1}{|t|}  | \langle e^{iHt}\varphi, A_{f'} e^{iHt}\psi \rangle | \leq \|\mathrm{ad}_{A_{f'}} H \| \, \|\varphi  \| \,\|\psi  \| \, .
\end{equation*}
Since $V\in C^1(X^2)$, we observe that for all $t\in {\mathbb R}$ and all $\varphi \in {\cal D}_X \subset {\cal D}(A_{f'}),$
\begin{equation*}
2\Re \langle i(\mathrm{ad}_X H) e^{iH t} \varphi, X e^{iH t} \varphi \rangle = 2\langle e^{iH t} \varphi, A_{f'} e^{iH t} \varphi \rangle + \bra e^{iH t} \varphi, i(\mathrm{ad}_{X^2}V) e^{iH t} \varphi \ket
\end{equation*}
where the last term on the RHS is uniformly bounded in $t$. The conclusion follows from Lemma \ref{propa1}. \ep

%******************************************
\subsection{Proof of Proposition \ref{propa3}}

As shown in the previous section, we observe that $H\in C^1(A_{f'})\cap C^1(X)$. We also deduce from Corollary \ref{Mourre-T2} and Lemma \ref{mourre-propa} that for all $\varphi \in {\cal H}$ such that $E_{\Lambda}(H)\varphi \in \mathcal{D}(A_{f'}) \cap \mathcal{H}_c(H)$
\begin{eqnarray*}
c_{\Lambda,f,|f'|^2}^{\sharp}\|E_{\Lambda}(H)\varphi\|^2 & \leq & \liminf_{t\rightarrow \pm \infty} t^{-1} \langle e^{iHt}E_{\Lambda}(H)\varphi, A_{f'} e^{iHt}E_{\Lambda}(H)\varphi \rangle \\
& \leq & \limsup_{t\rightarrow \pm \infty} t^{-1} \langle e^{iHt}E_{\Lambda}(H)\varphi, A_{f'} e^{iHt}E_{\Lambda}(H)\varphi \rangle \leq C_{\Lambda,f,|f'|^2}^{\flat}\|E_{\Lambda}(H)\varphi\|^2
\end{eqnarray*}
On the other hand, for all $t\in {\mathbb R}$ and all $\varphi \in {\cal H}$ such that $E_{\Lambda}(H)\varphi \in \mathcal{D}_X \cap \mathcal{H}_c(H) \subset \mathcal{D}(A_{f'}) \cap \mathcal{H}_c(H),$
\begin{eqnarray*}
2\Re \langle i(\mathrm{ad}_X H) e^{iHt}E_{\Lambda}(H) \varphi, Xe^{iHt}E_{\Lambda}(H) \varphi \rangle &=& 2 \langle e^{iHt}E_{\Lambda}(H)\varphi, A_{f'} e^{iHt}E_{\Lambda}(H)\varphi \rangle \\
&+& \langle i(\mathrm{ad}_{X^2} V)e^{iHt} E_{\Lambda}(H) \varphi, e^{iHt}E_{\Lambda}(H) \varphi \rangle
\end{eqnarray*}
since $V\in C^1(X^ 2)$. In particular, the last term on the RHS is uniformly bounded in $t$. The first statement follows from Lemma \ref{propa1}.

To prove the second statement, note first that for any smooth function $\Phi$ with compact support, $\Phi(H)\in C^1(A_{f'})\cap C^1(X)$. So, $\Phi(H)\varphi \in {\cal D}_X \cap \mathcal{H}_c(H)$ (resp. $\Phi(H)\varphi \in \mathcal{D}(A_{f'}) \cap \mathcal{H}_c(H)$) for all $\varphi \in {\cal D}_X \cap \mathcal{H}_c(H)$ (resp. $\varphi \in \mathcal{D}(A_{f'}) \cap \mathcal{H}_c(H)$). We deduce from Corollary \ref{Mourre-T3} and Lemma \ref{mourre-propa} that for all $\varphi \in \mathcal{D}(A_{f'}) \cap \mathcal{H}_c(H),$
\begin{eqnarray*}
c_{\Lambda,f,|f'|^2}\|\Phi(H)\varphi\|^2 & \leq & \liminf_{t\rightarrow \pm \infty} t^{-1} \langle e^{iHt}\Phi(H)\varphi, A_{f'} e^{iHt}\Phi(H)\varphi \rangle \\
& \leq & \limsup_{t\rightarrow \pm \infty} t^{-1} \langle e^{iHt}\Phi(H)\varphi, A_{f'} e^{iHt}\Phi(H)\varphi \rangle \leq C_{\Lambda,f,|f'|^2}\|\Phi(H)\varphi\|^2
\end{eqnarray*}
The conclusion follows as in the previous case substituting $E_{\Lambda}(H) \varphi$ by $\Phi(H)\varphi$, 
$C_{\Lambda,f,|f'|^2}^{\flat}$ by $C_{\Lambda,f,|f'|^2}$ and $c_{\Lambda,f,|f'|^2}^{\sharp}$ by $c_{\Lambda,f,|f'|^2}$. For the last statement, we refer to Corollaries \ref{Mourre-T3} and \ref{Mourre-T2}.

%***************************
\section{Complement on Laurent operators}

We can apply the techniques used previously to obtain similar results for (compact) perturbations of Laurent operators. We introduce them in this section. The details of the proofs are omitted. In the following $d\in {\mathbb N}$ is fixed.

\paragraph{Laurent operators.} Let $f\in L^{\infty}({\mathbb T}^d)$ and denote the sequence of its Fourier coefficients by $(\hat{f}_{\alpha})_{\alpha \in {\mathbb Z}^d}$. The Laurent operator $L_f$ associated to $f$ is the bounded discrete convolution operator defined by: $L_f : l^2({\mathbb Z}^d) \rightarrow l^2({\mathbb Z}^d),$
\begin{equation*}
L_f  \psi :=  \mathcal{F}(f)*\psi \, , \, \psi\in l^2({\mathbb Z}^d)\, .
\end{equation*}
In other words, for any $\beta \in {\mathbb Z}^d$, $(L_f \psi)_{\beta} =   \sum_{\alpha \in {\mathbb Z}^d} \hat{f}_{\alpha} \psi_{\beta-\alpha}$. Note that $\|L_f\|=\|f\|_{\infty}$ and that $L_f$ is unitarily equivalent to the multiplication operator by the function $f$ on $L^2({\mathbb T}^d)$: $L_f ={\cal F} f {\cal F}^*$. For any functions $f$ and $g$ in $L^{\infty}({\mathbb T}^d)$, any $c\in {\mathbb C}$: $L_{f+g}=L_f+L_g$, $L_{fg}=L_f L_g$, $L_{cf}=cL_f$, $L_f^*=L_{\bar{f}}$. In particular, $[L_f,L_g]=0$ and $L_1$ is the identity operator on $l^2({\mathbb Z}^d)$. If the function $f\in L^{\infty}({\mathbb T}^d)$ is real-valued, then the Laurent operator $L_f$ is self-adjoint and the spectral properties of $L_f$ are directly related to the properties of the function $f$. To name a few of them, we have that:
\begin{itemize}
\item the spectrum of $L_f$ is equal to the essential range of $f$. For example, if $f$ is continuous on ${\mathbb T}^d$, we have that $\sigma(L_f)=f({\mathbb T}^d)=$ Ran $f$, which is a connected and compact subset of ${\mathbb R}$.
\item $\lambda$ is an eigenvalue of $L_f$ if and only if $f^{-1}(\{\lambda\})$ has non zero Lebesgue measure
\item $L_f$ has purely absolutely continuous spectrum in a subset $\Lambda\subset {\mathbb  R}$
if and only if for any Borel set $N\subset \Lambda$ of zero Lebesgue measure, $f^{-1}(N)$ is also of measure zero.
\item $L_f$ has non-trivial singular continuous spectrum if and only if there exists a Borel set $N\subset {\mathbb R}$ of zero Lebesgue measure, such that $f^{-1}(N)$ has non-zero measure but $f^{-1}(\{\lambda\})$ is of zero measure for each $\lambda\in N$.
\end{itemize}
The reader will find a similar discussion in paragraph 7.1.4 \cite{abmg} for multiplication operators on ${\mathbb R}^d$.

\noindent{\bf Remark:} The shift operators $({\mathfrak S}_j)_{j\in \{1,\ldots,d\}}$ are examples of (unitary) Laurent operators on $l^2 ({\mathbb Z}^d)$: ${\mathfrak S}_j=L_{e^{i\theta_j}}$. Their action on the orthonormal basis of $l^2({\mathbb Z}^d)$ is given by: for all $\alpha \in {\mathbb Z}^d$, ${\mathfrak S}_j e_{\alpha} = e_{\alpha+\delta_j}$, where the vector $\delta_j \in {\mathbb Z}^d$ is defined by its coordinates: $(\delta_j)_n:= \delta_{jn}$, $n\in \{1,\ldots,d\}$. For all $(i,j)\in \{1,\ldots, d\}^2$: ${\mathfrak S}_j^*={\mathfrak S}_j^{-1}$ and $[{\mathfrak S}_i,{\mathfrak S}_j]=0$. For any $\alpha=(\alpha_1,\ldots,\alpha_d)\in {\mathbb Z}^d$, we write ${\mathfrak S}^{\alpha}={\mathfrak S}_1^{\alpha_1} \ldots {\mathfrak S}_d^{\alpha_d}$. If the symbol $f$ belongs to the
Wiener algebra ${\cal A}({\mathbb T}^d)$, then $L_f$ rewrites as a norm convergent series:
\begin{equation*}
L_f=\sum_{\alpha \in {\mathbb Z}^d} \hat{f}_{\alpha} {\mathfrak S}^{\alpha} \, .
\end{equation*}

We refer sometimes to the function $f$ as the symbol of the operator $L_f$. The sets of critical points of $f$ is
\begin{equation*}
\kappa_f = \{\theta \in {\mathbb T}^d; \text{f is not differentiable at } \theta \text{ or } \nabla f(\theta)=0\} \, .
\end{equation*}
In particular, if $f\in C^1({\mathbb T}^d)$ is real-valued, the set $\kappa_f$ is a compact subset of ${\mathbb T}^d$ and the set of thresholds $f(\kappa_f)$ is a compact subset of $\sigma (L_f)=\mathrm{Ran} f \subset {\mathbb R}$.

\paragraph{Conjugate operators.} Denote by $({\mathfrak X}_j)_{j=1}^d$ the family of linear operators defined on the canonical orthonormal basis of $l^2({\mathbb Z}^d)$ by: ${\mathfrak X}_j e_{\alpha} = \alpha_j e_{\alpha}$. The operators $({\mathfrak X}_j)_{j=1}^d$ are essentially self-adjoint on $\bra e_{\alpha}; \alpha \in {\mathbb Z}^d\ket$. We also denote by ${\mathfrak X}_j$, $j\in \{1,\ldots,d\}$ their respective self-adjoint extensions. In particular, for all $j\in \{1,\ldots, d\}$, ${\cal F}^*{\mathfrak X}_j {\cal F}= -i\partial_{\theta_j}$. This allows us to define the positive operator $\bra \mathfrak{X} \ket:= \sqrt{{\mathfrak X}_1^2+\ldots+{\mathfrak X}_d^2 +1}$ on
\begin{equation*}
{\cal D}_{\mathfrak{X}}=\cap_{j=1}^d {\cal D}({\mathfrak X}_j)=\{\varphi \in l^2({\mathbb Z^d}); \sum_{\alpha \in {\mathbb Z}^d}(1+\alpha_1^2+\ldots + \alpha_d^2) |\varphi_{\alpha}|^2 < \infty\}\, .
\end{equation*}
Using Fourier transform (\ref{fourier1}), we note that $L_f \in \cap_{j=1}^d C^1(\mathfrak{X}_j)$ if $f\in C^1({\mathbb T}^d)$. In this case, $L_f {\cal D}_{\mathfrak{X}} \subset {\cal D}_{\mathfrak{X}}$.

Let $g=(g_j)_{j\in \{1,\ldots,d\}}\subset C^2({\mathbb T}^d)$ a family of real-valued functions. We associate to $g$ a symmetric operator $\mathfrak{A}_g$ defined on ${\cal D}_{\mathfrak{X}}$ by:
\begin{equation*}
\mathfrak{A}_g := \frac{1}{2}\left( L_g \cdot \mathfrak{X} + \mathfrak{X} \cdot L_g \right)= \frac{1}{2} \sum_{j=1}^d L_{g_j} {\mathfrak X}_j + {\mathfrak X}_j L_{g_j}= \frac{1}{2} {\cal F} \left( g \cdot (-i\nabla) + (-i\nabla) \cdot g \right) {\cal F}^*\, .
\end{equation*}
Following the proof of Proposition 7.6.3 (a) in \cite{abmg}, we can show that the operator $i(g \cdot \nabla + \nabla \cdot g)$ is essentially self-adjoint on $C^2({\mathbb T}^d)$. We deduce that $\mathfrak{A}_g$ is essentially self-adjoint on ${\cal D}_{\mathfrak{X}}$. Its self-adjoint extension is also denoted $\mathfrak{A}_g$. Note that $\mathfrak{A}_g \bra \mathfrak{X} \ket^{-1}$ and $\mathfrak{A}_g^2 \bra \mathfrak{X} \ket^{-2}$ are bounded.

\begin{lem} Let $f\in C^3({\mathbb T}^d)$. Consider $g=(g_j)_{j\in \{1,\ldots,d\}}\subset C^2 ({\mathbb T}^d)$ a family of real-valued functions. Then $L_f\in C^1(\mathfrak{A}_g)$, $\mathrm{ad}_{\mathfrak{A}_g} L_f = L_{-ig\cdot \nabla f}$ and $\mathrm{ad}_{\mathfrak{A}_g} L_f \in C^2(\mathfrak{A}_g)$. In particular, $L_f \in C^3(\mathfrak{A}_{\nabla f})$ and $\mathrm{ad}_{\mathfrak{A}_{\nabla f}} L_f = -iL_{|\nabla f|^2}$.
\end{lem}
% $f\nabla\bar{f}=-\bar{f}\nabla f$
Given $f$, $g$ two real-valued continuous functions on ${\mathbb T}^d$ and any Borel set $\Lambda\subset {\mathbb R}$, $\Lambda \cap$ Ran $f \neq \emptyset$, let us define $c_{\Lambda,f,g}$, $C_{\Lambda,f,g}$, $c_{\Lambda,f,g}^{\sharp}$, $C_{\Lambda,f,g}^{\flat}$ as in Section 3.3 (with $f^{-1}(\Lambda)\subset {\mathbb T}^d$). Mourre inequality for the operator $L_f$ rewrites:
\begin{lem}\label{Mourre-L} Let $f\in C^3({\mathbb T}^d)$ be a non-constant real-valued symbol. Let $\Lambda \subset$ Ran $f$ be a Borel set. Then,
\begin{equation*}
C_{\Lambda,f,|\nabla f|^2}E_{\Lambda}(L_f) \geq E_{\Lambda}(L_f)(i\mathrm{ad}_{\mathfrak{A}_{\nabla f}} L_f) E_{\Lambda}(L_f) \geq c_{\Lambda,f,|\nabla f|^2}E_{\Lambda}(L_f)\, .
\end{equation*}
If $\overline{\Lambda}\subset$ Ran $f \setminus f(\kappa_f)$, then $c_{\Lambda,f,|\nabla f|^2}>0$.
\end{lem}
We deduce that if $V$ is a compact symmetric operator such that $V \in C^1(\mathfrak{A}_{\nabla f})$ and $\mathrm{ad}_{\mathfrak{A}_{\nabla f}} V$ is compact, then for $H=L_f+V$ and any real-valued $\Phi \in C^0({\mathbb R})$ vanishing outside the Borel set $\Lambda$, $C_{\Lambda,f,|\nabla f|^2}\Phi(H)^2 \gtrsim \Phi(H) (i\mathrm{ad}_{\mathfrak{A}_{\nabla f}} H) \Phi(H) \gtrsim c_{\Lambda,f,|\nabla f|^2}\Phi(H)^2$. Like Theorem \ref{toeplitz}, the next result is deduced by applying Mourre Theory (see Section 3.1):
\begin{thm}\label{laurent-glob} Consider a non-constant real-valued symbol $f\in C^3({\mathbb T}^d)$. Let $H=L_f+V$ with $V$ a compact symmetric operator defined on $l^2({\mathbb Z}^d)$ such that $V\in {\cal C}^{1,1}(\mathfrak{A}_{\nabla f})$. Then, $\sigma_{\text{ess}}(L_f+V)=\sigma_{\text{ess}}(L_f)=$ Ran $f$ and
\begin{itemize}
\item[(a)] given any Borel set $\Lambda$ such that $\overline{\Lambda}\subset$ Ran $f\setminus f(\kappa_f)$, $H$ has at most a finite number of eigenvalues in $\Lambda$. Each of these eigenvalues has finite multiplicity.
\item[(b)] a LAP holds for $H$ on Ran $f\setminus \sigma_{\text{pp}}(H)\cup f(\kappa_f)$ w.r.t $\mathfrak{A}_{\nabla f}$. $H$ has no singular continuous spectrum in Ran $f\setminus f(\kappa_f)$.
\end{itemize}
\end{thm}

\paragraph{An illustration.} The $d$-dimensional discrete Laplacian on ${\mathbb Z}^d$ described in \cite{sah} is actually the Laurent operator $L_f$ where $f$ is defined on ${\mathbb T}^d$ by $f(\theta)=2 \sum_{j=1}^d \cos \theta_j$: $L_f = \sum_{j=1}^d {\mathfrak S}_j + {\mathfrak S}_j^*$. We have that: $\sigma(L_f) =[-2d,2d]$, $\kappa_f=\{\theta\in {\mathbb T}^d; \forall j\in \{1,\ldots,d\}, \theta_j=0 \text{ or } \pi\}$ and $f(\kappa_f)=\{4k-2d; k\in \{0,\ldots,d\}\}$ is finite. If $V$ is compact symmetric and $H=L_f+V$, then $\sigma_{\text{ess}}(H)=\sigma_{\text{ess}}(L_f)=$ Ran $f$. We deduce from Theorem \ref{laurent-glob} that if in addition $V\in {\cal C}^{1,1}(\mathfrak{A}_{\nabla f})$, then
\begin{itemize}
\item Given any $\Lambda\subset {\mathbb R}$ such that $\overline{\Lambda}\subset [-2d,2d]\setminus f(\kappa_f)$, the operator $H$ has at most a finite number of eigenvalues in $\Lambda$. Each of these eigenvalues has finite multiplicity.
\item a LAP holds for $H$ on $[-2d,2d]\setminus \sigma_{\text{pp}}(H)\cup f(\kappa_f)$ w.r.t $\mathfrak{A}_{\nabla f}$. $H$ has no singular continuous spectrum.
\end{itemize}
We recover Theorem 2.1 in \cite{sah}.

\noindent{\bf Remark:} Note the existence of alternative methods for decaying potentials in random settings \cite{dss}, \cite{bourgain}.

In the context of Laurent operators, Proposition \ref{propa2} takes a slightly different form due to the commutation properties between $L_f$ and $L_{|\nabla f|^2}$ (see Lemma \ref{obs2}):
\begin{prop}\label{laurent-1} Consider a non-constant real-valued symbol $f\in C^3({\mathbb T}^d)$.
Then, for any $\varphi \in {\cal D}(\mathfrak{A}_{\nabla f}),$
\begin{equation*}
\lim_{t\rightarrow \pm \infty} \frac{1}{t} e^{-itL_f}\mathfrak{A}_{\nabla f} e^{itL_f}\varphi = L_{|\nabla f|^2}\varphi \, ,
\end{equation*}
and for any $\varphi \in {\cal D}_{\mathfrak{X}},$
\begin{equation*}
\lim_{t\rightarrow \pm\infty} \frac{1}{|t|} \|e^{itL_f} \varphi\|_{\mathfrak{X}} = \| L_{|\nabla f|}\varphi \| \, .
\end{equation*}
where $\|\varphi\|_{\mathfrak{X}} := \sqrt{\|\varphi\|^2+\sum_{j=1}^d\|{\mathfrak X}_j \varphi \|^2}$.
\end{prop}
The reader will also reformulate easily Propositions \ref{propa4} and \ref{propa3}:
\begin{prop}\label{laurent-2ub} Consider a non-constant real-valued symbol $f\in C^3({\mathbb T}^d)$.
Let $V$ be a bounded symmetric operator on $l^2({\mathbb Z}^d)$ such that $V\in C^{1}(\mathfrak{A}_{\nabla f})$ and $V\in \cap_{j=1}^d (C^1({\mathfrak X}_j)\cap C^1({\mathfrak X}_j^2))$. Let $H=L_f+V$. Then for any $\varphi \in {\cal D}_{\mathfrak{X}},$
\begin{equation*}
\limsup_{t\rightarrow \pm\infty} \frac{1}{|t|} \|e^{itH} \varphi\|_{\mathfrak{X}} \leq \sqrt{\|\mathrm{ad}_{\mathfrak{A}_{\nabla f}} H\|} \|\varphi \| \, .
\end{equation*}
\end{prop}

\begin{prop}\label{laurent-2lb} Consider a non-constant real-valued symbol $f\in C^3({\mathbb T}^d)$.
Let $V$ be a compact symmetric operator on $l^2({\mathbb Z}^d)$ such that $V\in C^{1}(\mathfrak{A}_{\nabla f})$, $V\in \cap_{j=1}^d (C^1({\mathfrak X}_j)\cap C^1({\mathfrak X}_j^2))$ and $\mathrm{ad}_{\mathfrak{A}_{\nabla f}} V$ is compact. Let $H=L_f+V$ and $\Lambda \subset$ Ran $f$ be a Borel set. Then,
\begin{itemize}
\item For any $\varphi \in {\cal H}$ such that $E_{\Lambda}(H)\varphi \in {\cal H}_c(H)\cap {\cal D}_{\mathfrak X},$
\begin{equation*}
\sqrt{c_{\Lambda,f,|\nabla f|^2}^{\sharp}} \|E_{\Lambda}(H) \varphi\| \leq \liminf_{t\rightarrow \pm\infty} \frac{1}{|t|} \|e^{itH}E_{\Lambda}(H) \varphi\|_{\mathfrak{X}} \leq \limsup_{t\rightarrow \pm\infty} \frac{1}{|t|} \|e^{itH}E_{\Lambda}(H) \varphi\|_{\mathfrak{X}} \leq \sqrt{C_{\Lambda,f,|\nabla f|^2}^{\flat}} \|E_{\Lambda}(H) \varphi\| \, .
\end{equation*}
\item For all $\varphi \in {\cal H}_c(H)\cap {\cal D}_{\mathfrak{X}}$ and all real-valued $\Phi\in C_0^{\infty}({\mathbb R})$ vanishing outside $\Lambda,$ 
\begin{equation*}
\sqrt{c_{\Lambda,f,|\nabla f|^2}} \|\Phi(H)\varphi \| \leq \liminf_{t \to \pm \infty} \frac{1}{|t|} \| e^{iHt}\Phi(H)\varphi \|_{\mathfrak{X}} \leq \limsup_{t \to \pm \infty} \frac{1}{|t|} \| e^{iHt}\Phi(H)\varphi \|_{\mathfrak{X}} \leq \sqrt{C_{\Lambda,f,|\nabla f|^2}} \|\Phi(H)\varphi \| \, .
\end{equation*}
\end{itemize}
If $\overline{\Lambda} \subset$ Ran $f\setminus f(\kappa_f)$, then $0<c_{\Lambda,f,|\nabla f|^2}^{\sharp}\leq c_{\Lambda,f,|\nabla f|^2}$.
\end{prop}

%*******************************

\end{document}